\documentclass[10pt, centertags]{amsart}

\newtheorem{thm}{Theorem}[section]
\newtheorem{prop}[thm]{Proposition}
\newtheorem{lem}[thm]{Lemma}
\newtheorem{cor}[thm]{Corollary}
\theoremstyle{definition}
\newtheorem{dfn}{Definition}
\theoremstyle{remark}

\newsavebox{\hjmtext}
\sbox{\hjmtext}{\small\textsc{Houston Journal of Mathematics}}
\newlength{\hjmlength}
\usepackage{amsfonts}
\usepackage{amssymb}
\usepackage{amscd}
\usepackage{amsmath}

\newcommand{\bdot}{  {\text{{\bf .}} } } 

\newcommand{\g}{\gamma}
\newcommand{\G}{\Gamma}

\newcommand{\nab}{\nabla}
\newcommand{\wh}[1]{\widehat{#1}}
\newcommand{\wt}[1]{\widetilde{#1}}

\newcommand{\FO}{{}^{\mathcal F}\mathcal O(M)}
\newcommand{\Dom}{Dom\'{\i}nguez}
\newcommand{\Ab}{{\mathcal A}_{\text{b}}}
\newcommand{\Cb}{C_{\text{b}}}
\newcommand{\Abp}{{\mathcal A}_{\text{b}}^\bot}
\newcommand{\Cbp}{C_{\text{b}}^\bot}
\newcommand{\kb}{\kappa_{\text{b}}}
\newcommand{\ko}{\kappa_{\text{o}}}
\newcommand{\phib}{\phi_{\text{b}}}
\newcommand{\phio}{\phi_{\text{o}}}

\newcommand{\gs}{ \text{{\bf gs}} }
\newcommand{\tF}{ {\widetilde {\mathcal F}}}
\newcommand{\tL}{ {\widetilde {\mathcal L}}}

\newcommand{\np}{\nab^\oplus}
\newcommand{\gp}{{}^\oplus\Gamma}
\newcommand{\ddt}{ {  {\text{d}\over \text{dt}}  } }

\newcommand{\tnu}{\widetilde{\nu}_f}

\let\leq=\leqslant
\let\geq=\geqslant
\let\pd=\partial

\let\isom=\approx
\newcommand{\visom}{\wr \wr}

\begin{document}
\title{An Application of Stochastic flows to Riemannian Foliations}
\author{Alan Mason}
\address{2316 Wright Circle\\
		 Round Rock, TX 78664
}
\email{Alanndg@aol.com}
\begin{abstract}
A stochastic flow is constructed on a frame bundle adapted to a
Riemannian foliation on a compact manifold.  The generator $A$ of
the resulting transition semigroup is shown to preserve the basic
functions and forms, and there is an essentially unique strictly
positive smooth function $\phi$ satisfying $A^*\phi = 0.$  This
function is used to perturb the metric, and an application of the ergodic theorem
shows that there exists a bundle-like metric for which the basic projection of
the mean curvature is basic-harmonic.
\end{abstract}

\maketitle
\theoremstyle{plain}
\newtheorem*{Prov}{Provision}

\section{INTRODUCTION}

To set the stage for the present work, let us begin by recalling the classical construction of Eells and Elworthy (see, e.g., [Bi], [IW]).  For $M$ a compact manifold with Riemannian metric $g$ and orthonormal frame bundle $\mathcal O(M)$, let $Y_i, i=1, \cdots, n = \text{dim }M$ be the
 canonical horizontal vector fields on $\mathcal O(M)$ for an affine connection $\nab$ that preserves the metric.  The idea of Eells and Elworthy is to consider the stochastic differential equation $dR_t = \sum_{i=1}^n Y_i(R_t) dw^i_t,\, R(0) = r_0$, and the associated
semigroup $(S_t f)(z) = \int_{\Omega}f\left(\pi(R(t, r, \omega))\right) P^W_0(d\omega)$. 
Here $f$ is a continuous function defined on $M$; the $w^i$ are the components of standard Brownian motion on $\Bbb R^n$; the differentials are taken in the Stratonovich sense; $z = \pi(r)$ is the projection of the initial frame $r$ at which the flow
 starts; and $P^W_0$ is Wiener measure on the space $\Omega$ of continuous paths starting at the origin. 
By using the globally defined vector fields $Y_i$ on the frame bundle and then projecting, this construction gives useful information while getting around the fact that $M$ is usually not parallelizable.  Thus, $S_t f(z)$ does not depend on the choice of frame $r$
over $z \in M$ and we have a well-defined object on $M$. 
Moreover, by Ito's formula and the Markov property of the flow, $S_t f$ satisfies the heat equation $\ddt S_t f(z) = A S_t f$ where the second-order elliptic operator $A$ depends on the connection $\nab$.

More precisely, the following facts are known.

	1) If $f$ is a smooth function on $M$ then $\wh A(f\circ\pi) = (A f)\circ \pi$, where $\pi: \mathcal O(M) \to M$ is the canonical projection.  Here $\wh A = {1\over2}\sum_{i=1}^n Y_i^2$, $A = {1\over2}\Delta + b$, $\Delta$ is the Laplacian for the Levi-Civit\`a connection, and $b$ is the so-called drift field.  Moreover, given any vector field $b$ on $M$, there exists a metric-preserving affine connection $\nab$ such that $b$ arises in this fashion.

	2) There exists a strictly positive smooth function $\phi$ on $M$ satisfying $A^*\,\phi = 0$, where $A^*$ is the $L^2$ adjoint of $A$.  Moreover, $\phi$ is unique up to multiplication by a constant (see Proposition 6.1 below).

	3) The (nonsymmetric) heat kernel for $A$ on functions is strictly positive and the ergodic theorem holds: $\lim_{t\to\infty} S_t f(z) = \int_M f d\mu$, where $\mu$ is the unique probability measure associated with $\phi$.

The above results hold quite generally, but as they stand there is no contact with geometry. It seems resonable that more information of a purely geometric nature should be obtainable if the drift field $b$ itself is geometrically well-motivated.

Pursuing this, let us suppose that we have some structure $\mathcal S$ on $M$, that is, a decomposition of $M$ into smooth orbits of a group action, say, or as leaves of a foliation, and suppose that this lifts to a structure $\wt{\mathcal S}$ on $\mathcal O(M)$ in the sense that each member $\wt{\Xi}$ of $\wt{\mathcal S}$ projects under $\pi$ to a member $\Xi$ of $\mathcal S$.  Finally, suppose that the metric-preserving connection $\nab$ used in the Eells--Elworthy
 construction has the following property: If $r_0$ and $r_1$ are two initial frames in $\wt{\Xi}_0$ then the corresponding flows $R(t, r_0, \omega)$ and $R(t, r_1, \omega)$ respect the structure after projection, i.e., for almost every $\omega$ and every $t\geq 0$, $\pi \circ R(t, r_0, \omega)$ and $\pi\circ R(t, r_1, \omega)$ belong to the same set $\Xi_t$.  Then the semigroup $S_t$ will also preserve the structure; for instance, if the function $f$ is constant on each set in the structure, the same will be true of $S_t f$.

Any connection $\nab$ with the above property will be of interest because the Eells--Elworthy construction then preserves useful geometric data.  In particular, the associated drift field $b$ will be a fundamental geometric object.  In the present paper we show that the above ideas can be fully implemented for Riemannian foliations, which seem almost tailor-made for our purposes.  Although the geometry has some novel features, the probabilistic techniques employed are standard.

We can now outline our results.  The key fact about Riemannian foliations that we need is the existence of bundle-like metrics.  These are used to lift the structure $\mathcal F$ to $\wt{\mathcal F}$ on $\mathcal O(M)$.  Lemma 2.1, a standard result for Riemannian submersions, depends essentially
on (3), while Lemma 2.2 uses nothing more than the characteristic property (2) of bundle-like metrics. The connection $\nab^\oplus = P\nab P + P^\bot\nab P^\bot$ is chosen for the Eells--Elworthy construction; that this is the right choice is shown in Lemmas 3.1 and 3.3.  The former lets us reduce to the adapted frame bundle $\FO$, while the latter, our main technical
result, shows that the transverse (deterministic) flows respect the foliation structure. 
This is not true for unrestricted flows, because Lemma 2.1 is valid only for tangent vectors $X$ perpendicular to the leaves. This means that the generator of the transverse transition semigroup arising from the construction will not be elliptic, but for the moment this causes no problems. 
\par
We next pass to the transverse stochastic flow in the standard way, and Lemma 3.4 shows that the associated semigroup preserves the basic functions.  We write $T_t$ for the transverse semigroup, reserving $S_t$ for the full semigroup.
Lemma 3.5 is a general result showing equality of transverse semigroups acting on basic functions under changes of metric.  All our results for $T_t$ are seen to hold already at the level of individual trajectories.  This is true in particular of Theorem 5.4, which is therefore merely a translation into heat-equation terms of the geometry of $\np$ using the Eells--Elworthy machinery.

Things become a little more interesting when we restrict our attention to functions and examine what facts 2) and 3) above have to say in our situation.  Here the need for ellipticity leads us to consider the full semigroup $S_t$.  Lemma 5.2 shows that $S_t f = T_t f$ for all basic functions $f$, even though the full flow does not respect the foliation structure; thus $S_t f = T_t f$ for basic $f$ does not follow by taking limits from a corresponding result that holds at the level of individual trajectories.  The proof of Lemma 5.2 uses the uniqueness of solutions of the heat equation and also, in an essential way, the fact that we have reduced to the subbundle $\FO$ and the fact that the transverse semigroup $T_t$ preserves the basic functions (Lemma 3.4).  
  
Lemma 5.1 establishes that the drift $b$ corresponding to $\np$ is just $\kappa/2$, where $\kappa$ is the mean curvature field; as expected, this is a fundamental geometric object.  We are now in a position to bring fact 3) to bear, leading to Theorem 6.2.  Remark 1 reflects the abundance of bundle-like metrics; for the purposes of Theorem 6.2, it would be sufficient just to dilate by $\phi$ or $\phib$.  Section 7 considers an example in some detail.  We remark here that since one cannot hope to actually calculate $\phi$ or its basic component $\phib$ explicitly, an essential role is played by the ergodic theorem as the only tool available for getting a handle on the behavior of $\phib$ when the bundle-like metric is varied.

\par
The author believes that an explicit geometric-probabilistic approach is the most natural, if not only, way to study the kinds of questions considered here in their full generality.  However, if one is just interested in Theorem 6.2, the question arises of whether the probability theory can be eliminated.  We will discuss this further at the end of the paper. 
\par
This work differs significantly from and largely supersedes the author's thesis [Ma], to which we can nonetheless refer for a few omitted proofs.
 
\section{THE ADAPTED FRAME BUNDLE AND ITS FOLIATION}

Let $M$ be a compact manifold of dimension $n$ equipped with a
foliation $\mathcal F$ of dimension $p$. There is an atlas of simple charts
$(U_\alpha, \phi_\alpha)$ on $M$ of the form
\[ \phi_\alpha:\;U_\alpha \isom \Bbb R^p\times\Bbb R^q \]
with distinguished coordinates
\[ \{z_j\} = \{x_i, y_{a-p}\}, \quad i=1, \cdots, p,\; a = p+1,
\cdots, n,\]
where the $x_i$ are along the foliation $\mathcal F$ and the $y_{a-p}$
are transverse to it.  Each subset $y = \text{const of } U$ is called a plaque and is contained in a leaf of $\mathcal F$; $q = \text{codim }\mathcal F$.

Let $q: z = (x, y)\mapsto \overline{z} := y$ also
denote the quotient map (defined locally on each chart), with differential
\[ q_*: T_zM \to \overline{Q}_z \equiv T_zM/T_z\mathcal F, \qquad X
\mapsto \overline{X}.  \]
Given a Riemannian metric $g$ on $M,$ we obtain a splitting
\[ TM = T\mathcal F \oplus Q \isom T\mathcal F\oplus \overline{Q} \]
of the exact sequence of bundles
\[ 0\to T\mathcal F \to TM \to \overline{Q}\to 0\]
where $Q = (T\mathcal F)^\bot,$ the orthogonal complement of $T\mathcal
F$ with respect to $g.$  If $(U', \phi')$ is another simple chart in the
atlas for $(M, \mathcal F)$, then the transition map $\phi'\circ\phi^{-1}$
on $U\cap U'$ is of the form
\begin{equation}(x, y) \mapsto (x'(x, y), y'(y)), \tag {1}
\end{equation} i.e., plaques go to plaques.
\par
We recall that a Riemannian foliation is one for which there exists an
atlas satisfying the following condition: the Jacobians ${\phi' \circ \phi^{-1}}_*$ define
maps $U \cap U'\to O(q),$ where $O(q)$ is the group of orthogonal matrices acting on $
\Bbb R^q.$  Equivalently, we can regard $\Bbb R^q$ as
a local model space equipped with a Riemannian metric $g_T$
which is preserved by the transition maps.  In general $g_T$ will not coincide
with the standard Euclidean metric on $\Bbb R^q$ and may have curvature; we will therefore
write $\overline{M/\mathcal F}$ rather than $\Bbb R^q$ for the local
model space.  A transverse covariant derivative $\nab^T$ on $\overline{M/\mathcal F}$
is uniquely determined by $g_T$ in the usual way by the Koszul
formula.  We will deal only with Riemannian foliations.
\bigbreak \par
\begin{dfn}  A vector field $\xi(z) =
\sum \xi_j(z) {\pd \over \pd z_j}$is said to be {\bf foliate} (or {\bf
projectable}) if it projects locally via $q$ to a vector field on
the local model space $\overline{M/\mathcal F},$ that is, if the
functions $\xi_j(z)$ for $j = p+1, \cdots, n$ depend only on the
$y$ coordinate in $z = (x, y).$
\end{dfn}

\begin{dfn} A form $\theta\in A^r(M)$ is said to be {\bf basic}
if for every $X\in T\mathcal F$ we have
\[ i_X(\theta) = 0\text{  and } i_X(d\theta) = 0,\]
where $i_X$ denotes contraction with $X.$ Thus $\theta$
is basic if and only if it involves only the transverse
coordinates $y$: $ \theta = \sum_K \theta_K(z) dz^K$
in terms of distinguished local coordinates $z = (x, y),$ where
$K = (k_1, \cdots, k_r)$ is an increasing multi-index with $k_1 >
p,$ and the coefficients $\theta_K$ depend only on $y.$  In
particular, a function is basic if and only if it is constant
along leaves.
We denote the spaces of basic functions and forms
by $C_{\text b}(M)$ and $\Ab(M),$ respectively.  The Riemannian metric $g$ defines
an $L^2$-projection $P_{\text{b}}$ onto the subcomplex of
basic forms and gives a decomposition $\theta = \theta_{\text{b}}
+\theta_{\text{o}}$ into basic and basic-orthogonal components.
\end{dfn}

\begin{dfn}The Riemannian metric $g$ on $M$ is {\bf bundle-like}
if and only if $\mathcal L_Z g = 0$ whenever $Z\in T\mathcal F$ is along the leaves; here $\mathcal L_Z$ denotes Lie derivative.  In other words,

\begin{equation}\tag{2} 
\begin{split}
&\text{ for any two local vector fields }X, Y \in (T\mathcal
F)^\bot, \text{ the function}\nonumber \\
& z\mapsto g_z(X, Y) \text{ is constant along the leaves wherever $X$ and $Y$ are foliate.} \\ 
\end{split}
\end{equation}
		
\end{dfn}

We will consider only bundle-like metrics $g$ that are compatible
with the given transverse metric $g_T$ in the following sense:
\begin{equation*} g_z(e, f) = (g_T)_{\overline{z}}(\overline{e},
\overline{f})\;\;\forall\;
     e, f \in T_z\mathcal F^\bot. \tag {3}
\end{equation*}
This is meaningful because the transverse metric $g_T$ is preserved under the
coordinate transformations in the defining atlas.
Such metrics can be constructed as follows. Given any Riemannian
metric $g'$ on $M,$ let $V \subset TM$ be the distribution defining the
foliation $\mathcal F,$ and let $P$ be the $g'$-orthogonal projection
on $V.$
Set $g(X, Y) = g'(PX, PY) + g_T(\overline{X}, \overline{Y})$ [Mo, Prop.\  3.3].
\par
There is an orthogonal splitting
\[TM = T\mathcal F \oplus T\mathcal F^\bot  \]
into vertical and horizontal subspaces. We write $P, P^\bot$
for the orthogonal projections on $T\mathcal F$ and $(T\mathcal
F)^\bot,$ respectively.  Because $g$ is compatible with $g_T$ (3), in each chart
$U_i$ with $q: U_i \to
\overline{M/\mathcal F}, \quad z \mapsto \overline{z} = y$ is a Riemannian
submersion onto the model quotient space, i.e., the
local quotient map $q$ gives an isometry $T_z\mathcal F^\bot \isom
T_{\overline{z}} \overline {M/\mathcal F}$.
\par
Passing to forms, we have a splitting $T^*M = T^*\mathcal F \oplus Q^*$ into
components along
and transverse to the leaves.  This induces a decomposition of the
$r$-forms on $M:$
\begin{equation}A^r(M) = \mathop\oplus_{u+v = r} A^u(Q)\otimes A^v(\mathcal F). \tag {4}
\end{equation}
There is a corresponding filtration, with forms in $A^{u, v} =
A^u(Q)\otimes A^v(\mathcal F)$ said to be of type $(u, v).$ With respect to
this filtration, the exterior derivative decomposes as $d = d_{1,0} + d_{0, 1}
+ d_{2, -1}$.
\par

Let $\mathcal O(M)\overset \pi \to M$ be the principal bundle of
orthonormal frames, and let $\FO$ be the subbundle of frames $r =
[z, (e_1, \cdots, e_p, e_{p+1}, \cdots, e_n)], \;\; z\in M,$
adapted to $\mathcal F.$  That is, the first $p$ vectors $e_i$ are
along the leaves, while the last $q$ are in $T\mathcal F^\bot.$
\par
In general, we say that a field of frames $r$ (i.e., a local
section of the bundle $\mathcal G\mathcal L(M)$ of all frames, or a
subbundle of it) is {\bf foliate} if each element $e_j$ is given
by a foliate vector field near $z.$ Expressing each $e_j$ as a
column vector in terms of the ${\pd\over\pd z_k},$ we see that a
frame in $\FO$ has the form
\begin{equation} r = \pmatrix A & B\\
         0 & C\\
     \endpmatrix.
\tag {5}\end{equation}  The $j$-th frame element is
\begin{equation} e_j = \sum^n_{k=1} e_j^k\,\pd/\pd z_k, \tag {6}\end{equation}
where $k$ labels the row and $j$ labels the column.
\par
Because the metric $g$ is bundle-like the Gram--Schmidt procedure, applied to
a preferred
basis \[\pd/\pd z_1, \cdots, \pd/\pd z_p, \pd/\pd z_{p+1}, \cdots, \pd/\pd z_n\]
in a simple chart, yields foliate frames, i.e., the elements $e_j\; (1\leq
j\leq n)$ are foliate. Gram--Schmidt thus creates a foliate local
orthonormal field of frames from a local chart.
\par

The following result will be needed in the construction of the flow. We
omit the straightforward proof, which uses (3) and the Koszul
formula for $\nab$ and ${\nab}^T$.
\begin{lem}If $X \in T_z\mathcal F^\bot$ then
\[ \overline{\left( P^\bot\nab_X P^\bot \pd/\pd z_l\right)_z }
= \left(
     \nab^T_{\overline{X}} \overline{ {\pd\over\pd
z_l}}\right)_{\overline{z}}.
\]
\end{lem}

\newcommand{\tX}{\widetilde{X}}

As the bundle-like metric $g$ varies, so do the spaces $\FO.$  We
will regard them as lying in $\mathcal G\mathcal L(M).$
\par
The adapted frame bundle $\FO \overset \pi \to M$ has a natural foliation
$\widetilde{\mathcal F},$ again of
dimension $p,$ which explicitly reflects the variation of the metric
$g$ along the leaves of $\mathcal F.$  The leaves of $\widetilde{\mathcal F}$
are of the form \[ \tL = \{r' = [z= (x, y), \vec{e'}] \;|\;
z \in\mathcal L, r' = \gs (r_0)\},\] where $\mathcal L$ is a leaf of $\mathcal F$
and $r_0 = [z_0=(x_0, y_0);
\vec{e}]$ is some reference frame based at a point $z_0 \in \mathcal
L.$ The components of $r' = \gs(r_0) = [z, \vec{e'}],\; z \in
\mathcal L,$ are by definition given by
\begin{equation}
\begin{align*}
e'_1 &= {e_1\over \|e_1\|_{g_z} }  \\
e'_2 &= {e_2-g_z(e_2, e'_1)e'_1\over
     \| e_2-g_z(e_2, e'_1)e'_1 \|_{g_z} }\\
     &\vdots\\
e'_{p+1} &=
{e_{p+1}-\sum^p_{j=1}g_z(e_{p+1}, e'_j)e'_j\over
     \| e_{p+1}-\sum^p_{j=1}g_z(e_{p+1},
     e'_j)e'_j \|_{g_z} } \\
     &\vdots\\
 \end{align*}
\tag {7}
\end{equation}

\newcommand{\rh}{\widehat{r}}
\newcommand{\eh}{\widehat{e}}

Here the reference frame $r_0$ is extended in the obvious way to
be a constant field in $\mathcal G\mathcal L(M)$ in a simple
chart about $z_0$: $r_0(z) = [z; \vec{e}],$ so that $e_j =
e^k_j(z_0)\pd_k$ is a constant vector field.  To make sense of this definition
 of $\tF$, we start with the fact that the Gram--Schmidt map $\gs$ is transitive: For $z, z', z''$ three points in a simple chart $U$, let
$ r'=\gs(r; z\to z'),\; r'' = \gs( r'; z'\to z''),\; \rh = \gs(r;
z\to z'')$; then $\rh = r''$. This leads to a global equivalence relation: $r\sim r'$ if and only if $r$ and $r'$ both lie over the same leaf $\mathcal L$ and there exist a chain of overlapping charts $U_i$ and frames $r_i\in \FO, z_i \equiv
\pi(r_i) \in \mathcal L \cap U_i, 0\leq i \leq N$, with $r = r_0, r' = r_N$, $z_i\in U_i\cap U_{i-1} \text{ for } 1\leq i \leq N$,
and $r_{i+1} = \gs(r_i, z_i\to z_{i+1})$ for all $i$. This equivalence
class of frames comprises the lifted leaf $\tL$ and defines the lifted foliation
$\tF$.  The transitivity of Gram--Schmidt ensures that there is no dependence on
the choice of reference frame $r_0 \in \tL$. It is easy to check that $\tF$ is a foliation,
and for each leaf $\tL$, $\pi: \tL\to \mathcal L$ is a covering map.
\begin{lem}The $C$ coordinates are constant along a leaf
$\tL$.\end{lem}
\begin{proof}  Since the $C$ coordinates of the first $p$ vectors are identically
zero for all frames $r$ in $\FO$, we start by considering
$e'_{p+1}$ in (7). Because $g$ is bundle-like and the local vector field $z\mapsto
e_{p+1}-\sum^p_{j=1}g_z(e_{p+1}, e'_j)e'_j$ is foliate and orthogonal to $T\mathcal
F,$ we have
\[
\| e_{p+1}-\sum^p_{j=1}g_z(e_{p+1}, e'_j)e'_j\| _{g_z} =
\|e_{p+1} - \sum^p_{j=1} g_{z_0}(e_{p+1}, e_j) e_j\|_{g_{z_0}} =
1.
\]
The assertion of the Lemma is now clear
for $e'_{p+1} = e_{p+1}-\sum^p_{j=1}g_z (e_{p+1}, e'_j)e'_j$.  Consider next
the numerator $e_{p+2}-\sum^{p+1}_{j=1}g_z(e_{p+2},
e'_j)e'_j$ of $e'_{p+2}.$  By (2), we have

\begin{align*}
g_z(e_{p+2}, e'_{p+1}) &=g_z\left(e_{p+2} - \sum^p_{k=1}
g_z(e_{p+2}, e'_k) e'_k,\quad e'_{p+1}\right)\\
     &= g_{z_0}(e_{p+2}, e_{p+1}) = 0.\\
\end{align*}

Thus $\|e_{p+2}-\sum^{p+1}_{j=1}g_z (e_{p+2}, e'_j)e'_j\| \equiv
1$ by the same argument used for $e'_{p+1},$ and hence $e'_{p+2}
= e_{p+2} - \sum^p_{j=1}g_z(e_{p+2}, e'_j)e'_j.$  Thus,
${e'}^k_{p+2} = e^k_{p+2}\text{ for all } k > p.$ Continuing in this way,
we obtain
${e'}^k_a = e^k_a \text{ for all }a, k > p.$
\end{proof}

Since the leaf $\tL$ is not globally contained in a simple chart, we need
to be more precise about the global meaning of Lemma 2.2.  To this end, let $C'$ be the
corresponding coordinates in an overlapping chart $U'$;
they are related to the coordinates $C$ by the Jacobian $\overline{J}(x, y)$
of the transformation (1), which is independent of the coordinates $x$
along the leaf $\mathcal L,$ given by $y = \text{const}.$ Since the leaf $\tL$
lies over $\mathcal L,$ we see that the $C'$ are constant along $\tL$ and given
by $C' = \overline{J}(x, y)\cdot C,$ for any value of $x$ corresponding to $z = (x, y), y = \text
{const},$ in the overlap $U\cap U'.$
Given two frames $r_0, r_1 \in \tL,$ we can join them by a path $\gamma$ in
$\tL$ and choose intermediate points $\rho_0 = r_0, \cdots, \rho_N = r_1$
on $\gamma$ such that the portion of $\gamma$ from $\rho_i$ to $\rho_{i+1}$
is contained in a simple chart $U_i,$ and
$\rho_i, \rho_{i+1}$ belong to the same plaque in $U_i.$  By following
along these plaques, we see how
the $C$ coordinates for $r_0$ are related to those for $r_1$ (in general,
there will of course be a dependence on the homotopy class of the path $\gamma$).
\par
On the other hand, by (7) the frame coordinates in $A$ transform by an
invertible matrix in $GL(p).$ The condition that the frames be orthonormal
at each point $z$ implies in particular:
\[g_z(A, B+C) = 0, \text{\qquad or\qquad} g_z(A, B) = -g_z(A, C)
\]
(in a convenient short-hand notation).  Thus $B$ is uniquely
determined by $C, \mathcal F,$ and the metric $g_z$; it does not
depend on $A,$ whose vectors merely span $T\mathcal F$.  As we
move along a leaf $\tL$, the metric varies and the $B$ components adjust
themselves so as to preserve orthogonality to $T\mathcal F$, the $C$
components remaining constant by Lemma 2.2.
\par
The structure group for $\FO$ is $G = O(p)\times O(q) \subset O(n)$.  A
frame $r = [z; \vec{e}]$ at $z\in M$ can be regarded as a map \[\Bbb
R^p\times\Bbb R^q \to T_zM,\;
(u, v) \mapsto \sum^p_1u_i e_i + \sum^n_{\alpha=p+1}v_{\alpha-p}e_\alpha.
\]
The action of $\g = \g'\times\g''$ is given by
\[(r\cdot\g)(u, v) = \sum^p_1(\g'\cdot u)_ie_i + \sum^n_{\alpha=p+1}
     (\g''\cdot v)_{\alpha-p}e_\alpha,
\]
where $(\g'\cdot u)_i = \sum^p_1 (\g')_{ij} u_j$ and so on. Thus,
the $j$-th frame element of $r\cdot\g$ is given by
\begin{equation}(r\cdot \g)_j
= \sum_i\g_{ij} e_i.\tag {8}\end{equation}

For $z_1, z_2 \in M$ and $r_1, r_2 \in \FO,$ we will write
\begin{equation} z_1 \sim z_2,\quad r_1 \sim r_2, \quad\text{ and } r_1\sim r_2
\text{ mod }O(p), \tag {9}\end{equation}
respectively, to mean that $z_1$ and $z_2$ lie on the same leaf
$\mathcal L$ of $\mathcal F;$ $r_1$ and $r_2$ lie on the same leaf $\tL$
of $\tF;$ and $r_2 \in \tL\cdot\gamma$ for some $\gamma \in
O(p),$ where $r_1 \in \tL.$  Clearly, $r_1\sim r_2 \text{ mod }O(p)$
implies $\pi(r_1) \sim \pi(r_2)$.
\par
Finally, for a given bundle-like metric $g$ on $M,$ we let $\nab$ denote
the Levi-Civit\`a connection on $M$ and set \[\nab^\oplus = P\nab P +
P^\bot\nab P^\bot. \]
Clearly, $\np$ preserves the metric $g$ since $\nab$ does.

\section{CONSTRUCTION OF THE FLOW}

To construct the flow we consider a simple chart $U$ with
coordinates $z = (x, y),$ in terms of which we have
\[ \np_{\pd_k} \pd_l = \sum_{i=1}^n \gp^i_{kl}\pd_i,\]
where $\pd_i = {\pd \over \pd z_i}$ and the $\gp^i_{kl}$ are the
Christoffel symbols.  Suppose that $i > p$ and $l
\leq p.$ Then $\np_{\pd_k}\pd_l = P\np_{\pd_k}\pd_l \in T\mathcal F,$
since $P^\bot \pd_l \equiv 0.$ Hence
\begin{equation} \gp^i_{kl} = 0 \quad \text{for } i> p, l \leq p.\tag {10} \end{equation}
Let $Y_a, \; 1\leq a\leq n,$ be the canonical horizontal vector
fields on $\mathcal G\mathcal L(M);$ they are uniquely determined by the
two conditions
\begin{equation*}
\begin{split}
 &i)\; Y_a \text { is horizontal for the connection }\np;\\
 &ii)\; \pi_*(Y_a |_r) = r(E_a)\in T_z(M)\\
\end{split}
\end{equation*}
for any frame $r\in \mathcal G\mathcal L(M), \pi(r) = z;$ here $E_a \in
\Bbb R^n$ is the canonical unit vector and we regard $r$ as a map
$\Bbb R^n \to T_z(M).$  We note that because $\np$ preserves the
metric, the $Y_a$ restrict to vector fields on the orthonormal
frame bundle $\mathcal O(M).$
\par
In terms of local coordinates $z, e^i_j$ on $\mathcal G\mathcal L(M)$ the standard
horizontal vector fields are given by [IW, Chap.\  V, Eq.\  (4.12)]
\begin{equation} Y_a = e_a^m\pd_m - \gp^i_{kl}\,e^k_a e^l_j \pd/\pd e^i_j; \tag {11} \end{equation}
all indices range from $1$ to $n,$ the ``vertical" coordinates
$e^i_j$ are given by $e_j = e_j^i\,\pd_i,$ and repeated indices
are summed.
\par

We fix a vector field $Y_a$ and consider the associated flow
${}_aR$ given by
\begin{equation*}\tag{12}
\begin{split}
     \ddt z^m(t) &= e_a^m(t) \\
     \ddt e^i_j(t) &= -\sum_{k,l}\gp^i_{kl}(z(t))e^k_a(t)e^l_j(t)\\
\end{split}
\end{equation*}
with initial condition ${}_aR(t=0) = r_0.$

\bigbreak \par
\begin{dfn} A flow $R(t, \cdot)$ will be said
to be {\it adapted}\; to $\mathcal F$ if $\pi\circ R(t, r_0)$
respects $\mathcal F$ in the following sense:
\[\pi \circ R(t, r_0) \text{ varies in a leaf } \mathcal L_t \text{
as } r_0 \text{ varies in } \tL.\]
This condition is weaker than requiring that the flow be foliate
for $\widetilde {\mathcal F}.$
We will say that $R(t, \cdot)$ is {\it weakly adapted}\; to $\mathcal F$ if:
\[\text{for every basic }f \in C_{\text b}(M), f(\pi(R(t, r_0))) \text{
is again basic,}
\]
for any choice of initial frame $r_0$ over $z \in \mathcal L.$  In
other words, given $z \in M,$ choose some frame $r_0\in \FO$ at
$z$ and let $r_0'$ vary in the leaf $\mathcal \tL$ containing $r_0;$
then $f(\pi(R(t, r_0')))$ is constant. 
\end{dfn}

In order for a flow $R(t, r_0)$ starting at $r_0 \in \FO$ to
be useful, it must preserve $\FO$ and be adapted to $\mathcal F$. The next two lemmas will show
that the flows ${}_aR, \,a=1, \cdots, n,$ have the necessary properties, even
though they are not foliate for $\widetilde{\mathcal F}.$
\bigbreak \par
\begin{lem} Let the flows ${}_aR,\; a = 1, \cdots, n,$ be as above.
Then each ${}_aR$ preserves $\FO$.
\end{lem}
\begin{proof}  Take $i > p, j\leq p,$ and pick $r_0 \in \FO$,
so that by (5), $e^i_j(t=0) = 0.$  We need to show that
$e^i_j(t) = 0$ for all $t.$  The right-hand side of the
second equation in (12) is zero at $t=0$ since $e^l_j(t=0) = 0$
unless $l \leq p,$ and by (10), $\gp^{i>p}_{k, l\leq p} \equiv
0.$ According to the theory of first-order differential
equations, if a flow starts at a point in a closed submanifold $N_1
\subset N$ and the vector field is tangent to $N_1$ at every
point in $N_1$, then the flow stays in $N_1;$ taking $N$ to be
$\mathcal G\mathcal L(M)$ and $N_1$ to be ${}^{\mathcal F}\mathcal G\mathcal L(M),$
the bundle of all frames with first $p$ vectors along $\mathcal F,$
we see that $e^{i>p}_{j\leq p}(t) = 0$ for all $t.$ Thus each
flow ${}_aR(t, \cdot)$ takes ${}^{\mathcal F}\mathcal G\mathcal L(M)$ to
itself.  Moreover, the vector fields $Y_a$ are horizontal for the
connection $\np,$ and (12) says precisely that each tangent
vector $e_j(t)$ is parallel along the curve $t\mapsto z(t).$
But parallel transport along $z(\cdot)$ preserves the
metric $g$ because $\np$ does; hence the ${}_aR$ also preserve
$\mathcal O(M).$  Therefore, they preserve $\FO = \mathcal O(M)\cap
{}^{\mathcal F}\mathcal G\mathcal L(M)$.
\end{proof}

The following immediate corollary deals with constant linear
combinations of the flows ${}_aR.$  The flow ${}_aR$ constructed
in Lemma 3.1 corresponds to the case $\vec{c} = E_a\in \Bbb R^n.$
\begin{cor} Consider the flow $R(t, \cdot,
\vec{c})$ given by the vector field $Y = \sum^n_1 c_i Y_i,$ where
the $c_i$ are constants.  Then $R$ preserves $\FO$.
\end{cor}

The next lemma is our main technical result.  Because Lemma 2.1 is not valid unless $X\in T_z{\mathcal F}^\bot$, we must limit ourselves here to transverse flows $R(t, \cdot,
\vec{c})$, those for which the first $p$ components $c_i, 1\leq i\leq p$, of $\vec{c}$ are zero.  

\begin{lem} Let $R(t, \cdot, \vec{c})$ be a transverse flow.  Then in the notation of (9), if $r_0 \sim r_1 \text{ {\rm mod }}O(p)$ we have
\[R(t, r_0, \vec{c}) \sim R(t, r_1, \vec{c}) \text{ {\rm mod }}O(p).
\]
In particular, $\pi(R(t, r_0, \vec{c})) \sim \pi(R(t, r_1,
\vec{c}))$, so $R$ is adapted to $\mathcal F$.
\end{lem}
\begin{proof} We give the proof in several steps, proceeding from local to global.
\par

1. The flow $R(t, \cdot, \vec{c})$ is defined by $Y = \sum_{a>p}^n c_aY_a$.  Thus
\[ Y = c_a e_a^m\pd_m - \gp^m_{kl}\,c_a e^k_a e^l_j \pd/\pd e^m_j,\]
where repeated indices are summed; $p+1\leq a\leq n, 1\leq m \leq n$, and so on.
Let us write
$X(t) = \sum_{a>p}^n c_a e_a(t)$, with $m$-th component $X^m(t) = \sum_{a>p}^n c_a e^m_a(t)$.  

According to (12), the equations for the flow in local coordinates read:

\begin{equation} \tag{*}
\begin{split}
\ddt e^{m}_{j}(t) &= -\gp^{m}_{k l}(z(t)) c_a e^k_a(t)e^l_j(t) \\
\ddt z^m(t) &= c_a e^m_a(t).\\
\end{split}
\end{equation}
We must show that $\pi\circ R(t, r_0)$ respects $\mathcal F.$ 
\par
Since $\pi_*$
kills the vertical directions and takes $\sum_{m=1}^p \sum_a c_a
e_a^m(t) \pd/\pd z_m$ to $T\mathcal F,$ we need only
check for each $m > p$ that $\sum_a c_a e^m_a(t, r_0) {\pd\over\pd z^m}$ is
foliate.  That is, for each $m, a > p$ there must be no dependence of $e^m_a(t,
r_0)$ on $r_0$ when $r_0$ varies locally along a leaf $\tL = \{r = [z, \vec{e}]\;|\; z\in \mathcal L, r = \gs(r_{\text{ref}})\}$ (by varying locally, we mean that $z_0 = \pi(r_0)$
remains within the chart $U$). 
\par
Thus we need to examine the above system of ordinary differential equations for $m>p$. Here our choice of the connection $\np$ is essential, as it
allows us to effectively decouple the coordinates in $C$ from
those in $A$ and $B.$  First of all, by (10) it follows that the terms on the right-hand side are zero unless $l > p$, and since all frames are in $\FO$, it follows that $j > p$ also, as otherwise $e^l_j(t) = 0$.  
In terms of the block decomposition in (5), the differential equations (*) for the
components in $C$ yield the transverse system of equations:
\begin{equation*} \tag{13}
\begin{split}
\ddt e^{m>p}_{j>p}(t) &= -\sum_{k>p, l>p, a>p}\gp^{m>p}_{k
l}(z(t))c_a e^k_a(t)e^l_j(t)\\
     &\qquad- \sum_{k \leq p, l > p, a>p}\gp^{m>p}_{k l}(z(t))c_a e^k_a(t)
e^l_j(t)\\
 &= -\sum_{l>p}\left(P^\bot \nab_{X(t)} P^\bot {\pd\over\pd
z_l}\right)^m e^l_j(t)\\
 &= -\sum_{l>p} \left( \nab^T_{\overline{X(t)}} \overline{
{\pd\over\pd z_l}}
     \right)^{m-p}_{\overline{z(t)}} e^l_j(t),\\
\ddt z^m(t) &= \sum_{a>p}^n c_a e^m_a(t).\\
\end{split}
\end{equation*}
In the first equality we have for emphasis separated out the terms with $k\leq p$; these correspond to the $B$ components of $X(t)= \ddt z(t)$. 
In the second equality we have used $m > p,$ so that
$(P\nab_{X(t)}P\,{\pd\over\pd z_l})^m = 0.$ The third equality
follows from Lemma 2.1 and involves only the coordinates
$\overline{z}, C.$  Thus the connection $\np$ has enabled us to split the $C$ coordinates off from the $A$ and $B$ coordinates.  By Lemma 2.2, the initial condition for
$\overline{z}, C$ remains the same as $r_0$ varies in $\tL.$
Hence the result follows since (13), taken for all $m>p$ and $j>p$, is a system of (nonlinear) first-order ordinary differential equations of the form $\ddt(\overline{z}(t), C(t)) =
F(\overline{z}(t), C(t)),$
where neither the initial condition nor $F$ depends
on the parameters $x$ along the leaf $\tL$. The solution $\overline{z}(t), C(t)$ is therefore independent of $r_0 \in \tL$ for all times $t$ provided the flow remains over $U$.

We conclude: Given frames $r_0, r_1 \in \tL$ with $z_0 = \pi(r_0), z_1=\pi(r_1)$ in $U$,
there exists $T > 0$ such that for all $t, 0\leq t \leq T$, we have
\[ C(R(t, r_0, \vec{c})) = C(R(t, r_1, \vec{c}))\]
and
\[\pi(R(t, r_0, \vec{c})) \sim \pi(R(t, r_1, \vec{c})).\]
By the definition of the lifted foliation $\tF,$ these two facts imply that
\[ R(t, r_0, \vec{c}) \sim R(t, r_1, \vec{c}) \text{ mod }O(p) \text{ for all } t, 0\leq t\leq T.\]
\par

We note that in addition to the transverse component which is well under
control, the flow $R$ also has vertical and longitudinal components about
which less can be said.  Because of the vertical component, even if $r_0$
and $r_1$ lie on the same leaf $\tL$, after a time $t$ we have only
$R(t, r_0, \vec{c}) \sim R(t, r_1, \vec{c}) \text{ mod } O(p)$; however, the vertical component
is of no consequence after we project by $\pi$. The longitudinal component, which for transverse flows is due to the bending of the leaves, on the other hand causes a drift along the leaves even after projection, and we must treat it together with the 
transverse motion in what follows.
\par
2. Suppose next that $r_0 \sim r_1\text{ mod }O(p)$ and $r_0$ lies on a leaf $\wt{\mathcal L}$;
then $r_1\cdot\g =:\widehat{r}_1 \in \tL$ for some $\g \in O(p).$ 
 Let $\tau$ be a 
path in $\tL$ joining $r_0$ and $\widehat{r}_1.$  We continue to work locally
and assume that the projection of $\tau$ under $\pi$ is contained in $U$.
By part 1),
\[ R(t, r_0, \vec{c}) \sim R(t, \widehat{r}_1, \vec{c})
\text { mod }O(p).\]
On the other hand, the system (12) now reads, with $Y_a$ replaced by $Y = \sum_{i=p+1}^n c_i  Y_i$:
\begin{equation*}
\begin{split}
 {dz\over dt} &= \sum c_i\,e_i,\cr
\np_{\dot z(t)}\vec{e}(z) &= 0,\cr
\end{split}
\end{equation*}
where $R(0) = r_0 = [z_0, \vec{e}_0]$ and $i = p+1, \cdots, n.$
Since for $h\in G = O(p)\times O(q)$ arbitrary we have $\sum_i (h^{-1}\vec{c})_i (\vec{e} h)_i = \sum_{i, j, k} h^{-1}_{ij}
c_j h_{ki} e_k = \sum_k c_k e_k,$ it is immediate from the form of
this equation that
\begin{equation}R(t, r\cdot h, \vec{c}) = R(t, r, h^{-1}\cdot\vec{c})\,h,
\tag {14} \end{equation}
where $h^{-1}\cdot\vec{c}$ denotes ordinary
multiplication of the vector $\vec{c}$ by the matrix $h^{-1}.$
This argument holds equally well for unrestricted $\vec{c} \in \Bbb R^n$ and
also establishes Eq.\  (19) below.  Taking $h = \g$, it follows that
\[ R(t, \widehat{r}_1, \vec{c}) = R(t, r_1\cdot\g, \vec{c})
= R(t, r_1, \g^{-1}\cdot\vec{c})\cdot\g.\]
Since $\g^{-1}\in O(p),$ we have $c_j = (\g^{-1}\cdot\vec{c})_j,
\;j = p+1, \cdots, n.$  Thus the transverse part (13) of the
system of equations is not changed by the action of $\g,$ so
\[ R(t, r_1, \g^{-1}\cdot\vec{c}) \sim R(t, r_1, \vec{c})\text{
mod }O(p)\]
is clear.
We conclude that there exists $T>0$ such that $R(t, r_0, \vec{c}) \sim R(t,
r_1, \vec{c}) \text{ mod }O(p)$ for all $t, 0\leq t\leq T$.
\bigbreak \par
3. Next let $r_0 \sim r_1\text{ mod }O(p),$ with no restriction that
$\pi(r_1)$ be in $U$.  As before, we have $r_1\cdot\g =: \widehat{r}_1 \in \tL$ for
some $\g \in O(p).$  Let $\tau$ be a path in $\tL$
joining $r_0$ and $\widehat{r}_1.$
We subdivide $\tau$ into segments, each of which projects under $\pi$ into
some simple chart, and apply step 2) to each segment.  We conclude that given
$r_0$ and $r_1$ with $r_0 \sim r_1\text{ mod }O(p)$, there exists $T>0$
such that
\[ R(t, r_0, \vec{c}) \sim R(t, r_1, \vec{c})\text{ mod }O(p)\]
for all $t, 0\leq t \leq T$.
\bigbreak \par

4. Finally, let $r_0, r_1\in \FO$ with $r_0 \sim r_1\text{ mod }O(p)$ be
arbitrary and define $T_0$ to be the supremum of all $t\geq 0$ such that
\begin{equation} R(t, r_0, \vec{c}) \sim R(t, r_1, \vec{c})\text{ mod }O(p).\tag {15} \end{equation}
We claim that $T_0 = \infty$.  If this is not so, then by the continuity of the
flow $R$ we may replace $t$ by $T_0$ in (15).  Applying part 3) to $R$ with
initial frames $r'_0 = R(T_0, r_0, \vec{c})$ and $r'_1 = R(T_0, r_1, \vec{c})$, and
using the group property of the flow: $R(t+s, r) = R(t, R(s, r))$, we see
that (15) holds for all $t$ between $0$ and some $T_1$ strictly greater than
$T_0$, contrary to the definition of $T_0$.
\end{proof}

Thus the transverse deterministic flows $ R(t, r, \vec{c})$ constructed above preserve
$\FO$ and are adapted to the foliation $\mathcal F$.  We next pass to the transverse stochastic
flow in the usual way by considering a dyadic decomposition $D_k, \;
k = 1, 2, \cdots,$ of the positive time axis into intervals $I_n = \{t\;|\;
n/2^k\leq t < (n+1)/2^k\},\; n = 0, 1, \cdots,$ and
imagining that the coefficients $c_i$ are randomly changed at times
of the form $t_n = n/2^k.$  By Lemma 3.3, the resulting flow
$R(t, \cdot)$, with the coefficients $c_i$ reshuffled in this
way, again preserves $\FO$ and is adapted to $\mathcal F.$  It is
possible to make sense of the limit as $k\to\infty,$ and the
result is called a stochastic flow.
\par
More precisely, consider the stochastic differential equation \begin{equation} dR_t = Y_i(
R_t) dw^i_t,\qquad R(0) = r_0,\tag {16} \end{equation}
where all differentials are
understood in the Stratonovich sense, and the $w^i, \;i=p+1, \cdots, n,$ are
the components of a
standard $q$-dimensional Brownian process $W$ on $\Bbb R^q.$  $W$
lives on $(\Omega_q, P_0^W),$ the space of all continuous paths $\omega:
[0, \infty] \to \Bbb R^q$
starting at $0,$ with the standard Wiener measure $P^W_0.$ It is known that
almost everywhere (with respect to $P^W_0$), each component $w^i$ is
H\"older continuous for any exponent $\alpha < 1/2,$ but is differentiable almost nowhere.
\par
There is a standard way to approximate the solution of (16) which involves
replacing the Stratonovich differentials in Eq.\  (16)
by a ``polygonal approximation'' on dyadic intervals:
\begin{equation} dR_t^{(k)} = \sum_{i=p+1}^n Y_i(R_t^{(k)}) \dot {w}^{i, k} \;dt,\quad
R^{(k)}(0) = r_0, \tag {17} \end{equation}
where
\[ \dot{w}^{i, k}(t) = 2^k\left(w^i(t^+_k) - w^i(t_k)\right),\]
with $t_k \equiv [2^kt]/2^k, t^+_k \equiv [1 + 2^k t]/2^k. $
These are ordinary differential equations on the frame bundle
with coefficients $c_i = \dot {w}^{i, k}$ constant on each dyadic
interval, and their integral curves define a flow of
diffeomorphisms.
\par
It is a fact that the sequence of maps $R^{(k)}(t, r_0, \omega)$
converges in probability to the solution $R(t, r_0, \omega)$ of Eq.\  (16),
uniformly on compact sets. Moreover, this convergence is
actually in the $C^m$ topology; hence there exists a subsequence
$R^{(k)}(t, r_0, \omega)$ of these diffeomorphisms which
converge, together with their derivatives with respect to $r_0,$
to the limit map $R(t, r_0, \omega),$ for almost every $\omega$ with respect to
$P_0^W.$ For this and related results, we refer to [Bi, Chap.\  1: Th.\ 2.1,
Th.\  4.1, and Th.1, p.\  71].
\par
It follows that the limit stochastic process $R_t$ will inherit any
properties of the approximating flows $R^{(k)}_t$ that persist under
closure. In particular, using Lemmas 3.1 and 3.3 the transverse stochastic flow (16)
constructed from the globally defined vector
fields $Y_i$ will be shown to preserve the adapted frame bundle $\FO$ and respect the
foliation $\mathcal F.$
\bigbreak \par

The flow (16) does not drop to a flow on $M,$ because of the
dependence on the choice of frame $r_0$ above $z_0 \in M.$
Nevertheless, the associated (transverse) transition semigroup $T_t$, defined on functions
$f \in C(M)$ by
\begin{equation} (T_t f)(z) = E[(f\circ\pi)(R(t, r, \cdot))] = \int_{\Omega_q}
f\left(\pi(R(t, r, \omega))\right) P^W_0(d\omega),\tag {18} \end{equation}
is independent of the choice of frame $r \in \FO$ over $z$.  This is because
the flow is equivariant:
\begin{equation} R(t, r\cdot\g; \omega) = R(t, r; \g^{-1}\cdot \omega)\cdot\g \tag {19} \end{equation}
cf. [IW, Chap.\  V, Eq.\  (5.7)].  Indeed, the transformation $\omega \mapsto 
\g\cdot \omega, \,(\g\cdot\omega)^i = \g^i_j \omega^j,$ leaves Wiener measure unchanged,
so that the probability law of the projection
$Z(t, z; \cdot) := \pi\circ R(t,r; \cdot)$
is independent of the choice of frame $r \in \FO$ above $z \in
M.$  Only this law, not the projected ``flow'' itself, is relevant in (18).
\bigbreak \par

\begin{lem}For almost every $\omega,$ the transverse stochastic flow $R(t,
\cdot, \omega)$ preserves $\FO$ and is adapted to the foliation $\mathcal F.$
In fact, there exists a $P_0^W$-negligible set $N$ such that for all $t
\geq 0$ and $\omega \notin N$ 
\begin{equation} R(t, r_0, \omega) \sim R(t, r_1, \omega) \text{ {\rm mod }}O(p) 
\text{ whenever } r_0 \sim r_1\text{ {\rm mod }}O(p). \tag {20} \end{equation}
\end{lem}
\begin{proof} We will need the case $m = 0$ of the following result [Bi, Th. 2.1]:\par
There exists a subsequence $n_k$ and a subset $N \subset \Omega$ with
$P_0^W(N) = 0$ such that for all $\omega \notin N,$
\[ R^{(n_k)}(t, \cdot, \omega) \text{ converges to } R(t, \cdot, \omega) \]
in the $C^m$ topology, uniformly on compact subsets of $\Bbb R^+ \times \FO.$
The approximations $R^{(k)}$ appearing here are the ones defined by (17). In
what follows we fix such a subsequence and for simplicity write $k$ for $n_k.$
That $\FO$ is preserved for all $\omega \notin N$ is clear, since each
approximation $R^{(k)}(t, \cdot, \omega)$ preserves $\FO$ and $\FO$
is closed in $\mathcal G \mathcal L(M).$
\par
Clearly, adaptedness is implied by (20), so it suffices to prove the
latter.  This follows from our previous results, which imply that the
approximations (17) satisfy (20).  Indeed, Lemma 3.3 applies and it is enough to
consider a composition $\Psi\circ\Phi$ of two diffeomorphisms, where
\[ \Phi = R(t, \cdot)
\text{ and } \Psi = R'(t', \cdot),\]
with $t = 1/2^k$ and $t'$ satisfying
$0\leq t'\leq 1/2^k.$ This composition corresponds to running (17) from time
zero to time $1/2^k + t',$ with initial point $r_0\in \FO;$ the flow $R'$
is obtained by reshuffling at time $t = 1/2^k$ the coefficients $c_i$ determining $R$,
as described after the proof of Lemma 3.3.
By Lemma 3.3 applied to $Y = \sum c_i Y_i,$ where the
$c_i$ are the constants for the flow $R$, we see that
$\Phi(r_0)\sim \Phi(r_1)\text{ mod }O(p).$
Now apply Lemma 3.3 again, this time to the reshuffled flow $R'$
with initial conditions $\Phi(r_0)$ and $\Phi(r_1)$, to conclude
that $\Psi(\Phi(r_0)) \sim \Psi(\Phi(r_1))\text{ mod }O(p)$ and
the approximating flows $R^{(k)}$ satisfy (20).  In particular,
$\pi(\Psi(\Phi(r_0))) \sim \pi(\Psi(\Phi(r_1))),$ so they are
adapted to $\mathcal F.$
\par
Finally, we need to show that the limit stochastic flow (16) on
$\FO$ satisfies (20).  This is not automatic, because the
leaves need not be closed.  Let $r_0 \sim r_1 \text{ mod }O(p)$ and repeat
the proof of Lemma 3.3, joining $r_0$ to $\widehat{r}_1$ by a path $\tau$
in $\tL.$  For fixed $t \geq 0$ and $\omega \notin N$ let us write $\Phi$
for the diffeomorphism $R(t, \cdot, \omega)$ of $\FO.$
Subdividing $\tau$ into small pieces and arguing on each
piece, we may suppose that $\tau$ is contained in a plaque in a
simple chart $\widetilde{U}$ and that the image of $\tau$
under $\pi\circ\Phi$ is contained in some simple chart $U$ with
distinguished coordinates $z = (x, y).$ Since the
$R^{(k)}(t, r, \omega)$ converge to $\Phi$ uniformly in $r \in \FO$ for all $\omega \notin N$,
for all sufficiently large $k$  
we have $\pi\circ R^{(k)}(t, r, \omega) \subset U$ for $r\in\tau$.
As shown in the previous paragraph, each $\pi\circ R^{(k)}(t, \cdot, \omega)$
takes plaques in $\FO$ to plaques in $M$, hence on taking limits we see that 
$\pi\circ\Phi(\tau)$  is contained in a plaque.  Moreover, the
$C$ coordinates of $R(t, r_0, \omega)$ and $R(t, \widehat{r}_1, \omega)$
coincide, since by the first part of this proof this is true for
the approximating flows $R^{(k)}(t, \cdot, \omega).$  From the
definition of $\tF$ (as in the proof of Lemma 3.3), it follows that
\begin{equation} R(t, r_0, \omega) \sim R(t, \widehat{r}_1, \omega)\text{ mod }O(p) \text{ for
almost every }\omega.
\tag {21} \end{equation}
\par

To finish, we observe that $r_1 = \widehat{r}_1\cdot \g$ for some
$\g \in O(p).$  Arguing as in the proof of Lemma 3.3, but using
Eq.\  (19) in place of (14), we obtain from (21) that $R(t, r_0,
\omega) \sim R(t, r_1, \omega) \text{ mod }O(p), \text{ a.e. }\omega.$
\end{proof}

In particular, $R(t, \cdot, \cdot)$ is weakly adapted to $\mathcal F,$ and hence
$T_tf$ given by (18) is basic whenever $f$ is.
\par
The next lemma establishes an important property of the transition
semigroup $T_t$ when $g$ is replaced by another
bundle-like metric $g'.$  We write $\FO$ and ${\FO}'$ for the
adapted orthonormal frame bundles for $g$ and $g',$ respectively;
the corresponding transverse transition semigroups are denoted by
$T_t$ and $T'_t.$  Recall that as remarked after Eq.\  (18), for
$f\in C(M)$, $ T_tf(z) = E[f(\pi R(t, r_0, \cdot))] $ and
$ T'_tf(z) = E[f(\pi R'(t, r'_0, \cdot))]$
do not depend on the choice of the initial frames $r_0 \in \FO$
and $r_0' \in {\FO}'$ over $z \in M.$
\bigbreak \par

\begin{lem}  For all $z\in M,$
we have
\begin{equation} T_tf(z) = T'_tf(z)\tag {22} \end{equation}
for all basic functions f. 
\end{lem}

\begin{proof}  By (18), (19), and the comment just before
Lemma 3.4, we may replace the initial frame $r_0'\in {\FO}'$ by
$r'_0\cdot\gamma,\,\g \in G = O(p)\times O(q).$  By (3), we can
choose $\g\in O(q)$ so that, in the notation of (5), the frame
coordinates $C_0'$ for $r_0'\cdot\g$ coincide with $C_0$ for
$r_0.$
\par

We begin by arguing locally within a coordinate chart $U_1$. Recalling (13),
we get the transverse systems of differential equations for the two transverse deterministic flows $R$ and $R'$ in local coordinates:

\begin{equation*}\tag {23}
\begin{split}
\ddt e^{m>p}_{j>p}(t) &= -\sum_{l>p, a>p, k}\gp^{m>p}_{k
l}(z(t))c_a e^k_a(t)e^l_j(t)\cr
 &= -\sum_{l>p}\left(P^\bot \nab_{X(t)} P^\bot {\pd\over\pd
z_l}\right)^m e^l_j(t),\cr
\ddt z^m(t) &= X^m(t)
\end{split}
\end{equation*} 

and
\begin{equation*} \tag {24}
\begin{split}
\ddt {e'}^{m>p}_{j>p}(t) &= -\sum_{l>p, a>p, k}{\gp'}^{m>p}_{k
l}(z'(t))c_a {e'}^k_a(t){e'}^l_j(t)   \cr
 &= -\sum_{l>p}\left({P'}^\bot \nab_{X'(t)} {P'}^\bot {\pd\over\pd
z_l}\right)^m {e'}^l_j(t),\cr
\ddt {z'}^m(t) &= {X'}^m(t)
\end{split}
\end{equation*}
In writing (24) we use the direct-sum connection ${\np}'$ for the
metric $g'$ on $M$ and the associated canonical vector fields $Y'_i$; ${P'}^\bot$ is 
the orthogonal projection on $(T\mathcal F)^\bot$ for $g'.$   Recall that $X(t) = \sum_{a>p}^n
c_a e_a(t)$, and we define similarly $X'(t) = \sum_{a>p}^n c_a e'_a(t)$.  

By Lemma 2.1, we have (as in the first part of the proof of Lemma 3.3)
\begin{equation*}
\begin{split}
 \overline{
          P^\bot\nab_{X(t)} P^\bot {\pd\over\pd z_{l>p}}
         }
     &= \nab^T_{ \overline{X(t)}} \overline{ {\pd\over\pd
z_l}}\quad\text{(at $\overline{z}(t)$)}\cr
     \overline{
          {P'}^\bot\nab'_{X'(t)}{P'}^\bot {\pd\over\pd
          z_{l>p}}
     }
     &= \nab^T_{ \overline{X'(t)}} \overline{ {\pd\over\pd
z_l}}\quad\text{(at $\overline{z'}(t)$)},\cr
\end{split}
\end{equation*}
where $\nab^T$ denotes the Levi-Civit\`a connection for the
transverse metric $g_T$ on the local model space
$\overline{M/\mathcal F}$.
Thus the form of the two equations (23), (24) for the coordinates
$(\overline{z}, C)$ and $(\overline{z'}, C')$ is identical; since the initial
conditions coincide, we see that $(\overline{z(t)}, C(t)) = (\overline{z'(t)}, C'(t))$.
\par
Next, we must globalize this result.  The difficulty is that although the
transverse parts of $g$ and $g'$ are the ``same'' by (3), there is no
correlation in the variation of the longitudinal parts of $g$ and $g'$
as we move along a leaf. This results in a longitudinal drift of the two
flows relative to one another which must be treated here.
\par
Fix some time $t > 0$ such that for all $ 0 \leq \tau \leq t$,
both $\pi(R(\tau))$
and $\pi(R'(\tau))$ lie within the chart $U_1,$ while $\pi(R'(t))$
also lies in an overlapping chart $U_2.$  The initial frames for $R, R'$
are $r_0 \in \FO$ and $r'_0 \in {\FO}'.$
Before starting up the flows, we were free to replace $r'_0$ by
$r'_0\cdot\gamma,\;\gamma\in O(q),$ so that its initial $C$ coordinates
$C'$ agreed with those of $r_0.$  As the flows evolve in time, however,
it is essential that we not do this again as this would change the transverse
equations (24) for $R'(t),$  which is not allowed.
\par
By the part of Lemma 3.5 already proved, we have
\begin{equation}C'(t) = C(t) \tag {25} \end{equation}
using the coordinates in the chart $U_1$, and the projections
$z_t=\pi(R_t)$ and $z'_t=\pi(R'_t)$ lie on the same leaf $\mathcal L_t$ of
$\mathcal F$.  (Here we write $R_t$ for $R(t, r_0)$ and similarly for $R'_t$.)
Let $\sigma$ be a path in $\mathcal L_t \cap U_1$ from $z_t$ to $z'_t$ and
let $\widetilde{\sigma}$ be the lift of $\sigma$ starting at $R_t$ and
contained in $\tL_t$.  The endpoint $A_t$ of $\widetilde{\sigma}$
satisfies $\pi(A_t) = z'_t = \pi(R'_t).$  Let ${}^{\text{tr}}R: s\mapsto
R(s, A_t)$ denote the ``translated'' flow with initial value $A_t, 0 \leq s$.
By Lemma 2.2 applied to the metric $g$, bundle $\FO$, and lifted foliation $\tF$,
\[C(A_t) = C(t)\]
because $\sigma $ lies within the chart $U_1$.   Thus, by Eq. (25) we have
\begin{equation} C(A_t)= C'(t) \tag {26} \end{equation}
in terms of the coordinates for the chart $U_1,$ and therefore also in
terms of the coordinates in the overlapping chart $U_2$ (recall the
discussion after Lemma 2.2).
\par
The essential point is that by Eq. (26), the new initial points
$R'_t$ and $A_t$  are already ``in register'' in terms of the coordinates
of chart $U_2,$ so no further application of $\gamma\in O(q)$ is necessary.
Letting the flows develop from $A_t= {}^{\text{tr}}R(s=0) $ and
$R'(0, R'_t)$ for a time $s > 0$ small enough so that we remain in $U_2$,
we obtain (using the semigroup property of the flows and the notation of (9)):
\[ \pi (R_{t+s}) \sim \pi ({}^{\text{tr}}R_{s}) \sim \pi (R'_{t+s}). \]
The first relation holds by Lemma 3.3 applied to $R$, and the second follows
by an another application of the first part of the proof of Lemma 3.5,
this time within the chart $U_2.$
\par
Thus we can use Lemma 3.3 to translate the flow $R_t$ along $\tF$, compare
the translated flow with $R'_t$ in some other chart, and deduce that
$ \pi(R_t) \sim \pi(R'_t)$ for all times $t\geq 0.$
\par 
The next step is to treat the approximating flows $R^{(k)}(t,
\cdot)$ in (17), which is done by considering composites of flows
corresponding to vector fields $Y = \sum c_i Y_i$ with initial
conditions $r_0 \in \tL\cdot O(p).$  The argument is the same as
in the proof of Lemma 3.4.
\par
Thus the approximating flows satisfy $\pi(R^{(k)}_t) \sim \pi({R_t'}^{(k)})$ for all $t \geq 0$, and the analogous result for the stochastic flows holds for almost every $\omega$ on passing to the limit. The equality (22) now follows from (18).
\end{proof}
\section{EXTENSION TO FORMS}

Let $u$ be a tensor of type $(a, b).$  In terms of the local
coordinates $z_1, \cdots, z_n,$ $u(z)$ is given in terms of its
components $u(z)^K_L$ by  \[u(z) = u(z)^K_L \pd_K\otimes dz^L,\]
where $K = (k_1, \cdots, k_a)$ and $L = (l_1, \cdots, l_{\text {b}})$ are
multi-indices of degree $a$ and $b;$ $\pd_K \equiv {\pd\over\pd
z_{k_1}}\otimes\cdots\otimes {\pd\over\pd z_{k_a}}$ and
$dz^L\equiv dz^{l_1}\otimes\cdots \otimes dz^{l_b}.$
\par
In terms of frames $r = [z; \vec{e}]$ we can write
\begin{equation} u(z) = F^I_{uJ}(r)\,e_I\otimes e_*^J = F^I_{uJ}(r)\,e^K_I\,f^J_L\;\pd_K\otimes dz^L, \tag {27}\end{equation}
where $I, J$ are multi-indices, and $e_I \equiv e_{i_1} \otimes
\cdots\otimes e_{i_a},$ and so on.  The coordinates $e^i_k,
f^k_i$ of the $k$-th frame vector $e_k$ and the $k$-th vector
$e^k_*$ of the dual frame are defined by
\begin{equation}  e_k = e_k^i {\pd\over\pd z_i}, \qquad e_*^k = f^k_i dz^i;
\tag {28}\end{equation}
the matrix $(f_i^j)$ is the inverse of $(e_i^j).$  If $r = [z;
\vec{e}]$ is expressed in block form as in Eq.\  (5), then
\[ (e^i_j) = \pmatrix A & B\\
          0 & C\\ \endpmatrix
\qquad\text{ and\qquad}
  (f^i_j) = \pmatrix A^{-1} & -A^{-1}BC^{-1}\\
          0 & C^{-1}\\ \endpmatrix
\]
The functions $F^I_{uJ}$ are well-defined on the entire frame
bundle; however, the components $e^K_I, f^J_L$ in (27) are
defined only with reference to the local chart $\{z_j\}.$  Observe that
the definition (28) for $e^k_*$ involves the transpose of $(f^i_j)$; thus
we regard $e_k$ as the $k^{\text{th}}$ column vector of $(e^i_j)$ and
$e^k_*$ as the $k^{\text{th}}$ row vector of $(f^i_j)$.  The $e_k$ with
$1\leq k \leq p$ span $T\mathcal F = \text{span}\{\pd/\pd z_i\}, 1\leq i \leq p$,
while the $e^k_*$ with $p+1\leq k \leq n$ span the transverse space
$Q^* = \text{span}\{ dz^a\}, p+1\leq a\leq n$.
\par
The collection of functions $\{F^I_{uJ}\}$ on the frame bundle is called
the {\it scalarization}\;of $u$ and is equivariant (see, e.g., [IW, p.\  280]
or [BGV, p.\  24]).  That is, \begin{equation}F^\bdot_{u\bdot}(r\cdot \g) = F^\bdot_{u\bdot}
(r)\cdot (\g^\otimes)^{-1}, \tag {29}\end{equation} where $r\cdot\gamma$ is given by (8).
\par

Conversely, if (29) holds for some collection $\{F^I_J\}$ of
functions, then there exists a unique tensor $u$ of which
$\{F^I_J\}$ is the scalarization. We have
\begin{equation*} \tag {30}
\begin{split}
u(z)^K_L &= F^I_{uJ}(r) e^K_I f^J_L,\cr
F^I_{uJ}(r) &= u(z)^K_L e^L_J f^I_K. \cr
\end{split}
\end{equation*}

We now specialize to the case when $u = \theta(z) =
\theta(z)_J\,dz^J$ is an $m$-form and consider only frames $r
\in \FO.$
\bigbreak \par

\begin{lem}  \it $\theta$ is basic if and only if:
\par
\qquad $i)\quad\text{each } F_{\theta J} \text{ is constant along
}
\tL\cdot O(p)
     \text{ ($\tL$ a leaf of }\tF) $
and
\par
\qquad $ii)\quad F_{\theta J}(r) = 0 \text{ whenever any index }
j_\nu\leq p.$
\par\noindent 

In other words, $\theta$ is basic
if and only if the $F_{\theta
J}$ depend only on the $C$ coordinates for $J > p$ and vanish
otherwise.
\end{lem}

\begin{proof}  The straightforward proof [Ma] is based on Lemma 2.2.
\end{proof}

\bigbreak \par

Given a form $\theta$ with scalarization $\{F_{\theta J}\},$ we set
 \begin{equation} U_J(t,
 r_0) = E[F_{\theta J}(R(t, r_0, \omega))] \equiv \int_{\Omega_q} F_{\theta J}(R(t,
 r_0, \omega)) P^W_0(d\omega). \tag {31}
 \end{equation}
By (19), the transverse flow $R$ is $G = O(p)\times
O(q)$-equivariant.  Since
$\{F_{\theta J}(\cdot)\}$ is equivariant (29), the same is true of
$\{U_J(t, \cdot)\}$ for each $t \geq 0,$ because $\omega\mapsto \gamma\cdot
\omega$ leaves the measure $P_0^W$ unchanged. By the observation made after
(29), it follows that there exists a unique $m$-form $\theta(t, z_0)$ of
which $\{U_J(t, r_0)\}$ is the scalarization.
The action of the transverse semigroup $T_t$ on forms is defined by
\begin{equation} (T_t\theta)(z) = \theta(t, z). \tag {32}\end{equation}
We have
\begin{lem}$T_t\theta$ is basic whenever
$\theta$ is.
\end{lem}

\begin{proof} This follows from Lemmas 3.4 and 4.1. 
\end{proof}

	We note here that the extension (32) of $T_t$ to differential forms is
easily seen to preserve the filtration (4).

\section{THE HEAT EQUATION}
	
We now consider, in addition to the transverse semigroup $T_t$ constructed above, the full
semigroup $S_t$ constructed as in (18), but using the full stochastic flow $R(t, r, \omega)$ constructed as described after Lemma 3.3 from the unrestricted deterministic flows $R(t, r, \vec{c})$, for which $\vec{c} \in \Bbb R^n$ is arbitrary; thus in (18), $\Omega_q$ is replaced by $\Omega_n$.  The infinitesimal generator of $S$ is elliptic, as required for strict positivity of the heat kernel and ergodicity, which we need in Section 6.  However, because the full flow does not respect the foliation, it is not clear that $S_t$ preserves the basic functions, though this crucial property holds for $T_t$ (Lemma 3.4).  Nevertheless, it is a remarkable fact that after the averaging over $n$-dimensional Wiener measure is performed to get $S$ we have $S_t f = T_t f$ for all basic functions $f$.  In the present section we prove this result and examine some properties of the infinitesimal generators.

We begin by recalling the fundamental result [IW, Chap.\  V, Th.\  3.1] that the transition
semigroups $T_t$ and $S_t$ defined by (18) give solutions to the heat
equation. Namely, set $\tnu(t, r) = S_t f(t, r) \equiv
E[f(R(t, r, \cdot))]$ for any $f \in C^\infty(\FO)$; then $\tnu$
satisfies the partial differential equation
\begin{equation}{\pd \tnu\over\pd t} = {1\over2}\sum^n_1 Y^2_k\,\tnu, \qquad
\tnu(0, r) = f(r).
\tag {33}\end{equation}
Let us write \begin{equation}\widehat{A} \equiv {1\over2}\sum^n_1 Y^2_k. \tag {34}\end{equation}
In the corresponding equation for the transverse semigroup $T_t$, $\wh{A}$ is replaced by $\wh{A^\bot}$, the summation over $k$ now going from $p+1$ to $n$.

The proof of the next lemma is an application of [IW, Chap.\  V,
Eq.\  (4.33)]; indeed, Ikeda and Watanabe show that any drift vector field
$\vec{b}$ on $M$ can be obtained by using a suitable affine
connection $\nab$ on $M$ that preserves the metric but has
nonzero torsion in general [IW, Prop.\  V.4.3]. The direct sum
connection $\np$ used here preserves the metric, and we will now see that its
torsion is such that the drift field $\vec{b}$ is just ${1\over2}\kappa,$
where $\kappa$ is the mean curvature field.
\begin{lem} For $f\in C^\infty(M),$ consider
the lift $f\circ\pi$ to $\FO,$ and let $\widehat{A}$ be as in
(34).  Then
\begin{equation}\widehat{A} (f\circ\pi) = (A f)\circ\pi, \tag {35}\end{equation}
where
\begin{equation} A = {1\over2}\Delta_M + {1\over2}\kappa. \tag {36}\end{equation}
Here $\Delta_M = -\delta d = +g^{ij}{\pd\over\pd z_i}{\pd\over\pd
z_j} - g^{ij}\G^k_{ij}{\pd\over\pd z_k}$ is the Laplacian for the
given bundle-like metric $g.$
\end{lem}

\begin{proof}  The drift field $\vec{b}$ is given in local coordinates by
\begin{equation} b^i = {1\over2}g^{km}\left( \Gamma^i_{km} -
{}^\oplus\Gamma^i_{km}\right), \tag {37}\end{equation}
where $\Gamma^i_{km}$ and ${}^\oplus\Gamma^i_{km}$ are the
Christoffel components for the Riemannian and direct-sum
connections, respectively.  Moreover, (35) holds with $A =
{1\over2}\Delta_M +\vec{b}$, see [IW, Chap.\  V, Eq.\  (4.33)].
\par

To show (36), pick $z\in M$ and a simple neighborhood $U \owns z$
in $M$ with coordinates $z_a,$ such that the $z_a = x_a$ with
$1\leq a\leq p$ are along $\mathcal F$ while the $z_b = y_{b-p},
\,p+1\leq b \leq n,$ are transverse. By definition, the mean
curvature is the vector field given by
\begin{equation} \kappa = \sum^p_{a=1}\sum^n_{b=p+1}\;g(\nab_{e_a}\,e_a, \;
e_b)\,e_b,\tag {38}\end{equation}
for any local orthonormal frame $\{e_i\}$ with $e_a$ in $T\mathcal F$
and $e_b$ in $(T\mathcal F)^\bot.$
We will take the $e_i,\;1\leq i \leq n,$ to be obtained by applying the
Gram--Schmidt procedure to \[\pd/\pd z_1, \cdots, \pd/\pd z_p,
\pd/\pd z_{p+1}, \cdots, \pd/\pd z_n,\] in the given order.  We have seen
that because the metric $g$ is bundle-like, the $e_i$ are foliate
(recall the discussion preceding Lemma 2.1). Since the vector field
$\vec{b}$ is tensorial, in (37) we can work with the local field of
orthonormal frames $\{e_i\}$ just constructed and obtain
\begin{equation*} \tag {39}
\begin{split}
 2b^i &= \sum_k (\nab_{e_k}e_k - {}^\oplus\nab_{e_k}e_k,\; e_i)\cr
    &=\sum_{k\leq p}(e_i,\;P^\bot\nab_{e_k}e_k) +
\sum_{k>p}(e_i,\; P\,\nab_{e_k} e_k).\cr
\end{split}
\end{equation*}

\bigbreak \par
We consider the two cases $i>p$ and $i\leq p$ separately.

For $i>p$ we have $2b^i = \sum_{k\leq p}g(e_i,\,\nab_{e_k}e_k) =
\kappa^i$ by (38).
\par

For $i\leq p,$ (39) reduces to
\[ 2b^i = \sum_{k>p}g(e_i, \,\nab_{e_k} e_k).\]
By the Koszul formula,
\[ 2g(\nab_{e_k}e_k,\; e_i) = 2 g(e_k, [e_i, e_k]),\]
which is zero because $e_{k>p}$ is foliate, i.e., $[e_{i\leq p},
e_k] \in T\mathcal F.$
We conclude that $\vec{b} = {1\over2}\kappa$.
\end{proof}
\bigbreak \par

For $f\in C^\infty(M)$ and $z\in M,$ let us write
$\nu_f(t, z) \equiv \widetilde {\nu}_{f\circ\pi}(t, r) =
E[f\circ\pi(R(t, r, \cdot))],$ where $\pi(r) = z$ and we are using the full flow $R$;
by the discussion after (18) this is well-defined, i.e., independent of
the choice of frame $r$ over $z.$  Since
$\widetilde{\nu}_{f\circ\pi}(t, r) = \nu_f(t, \pi(r)),$ it
follows from equation (33), with $f$ replaced by $f\circ\pi,$ and
the relation (35): $\widehat{A}(\nu_f\circ\pi) = (A\nu_f)
\circ\pi$, that $\nu_f(t, z) = (S_t f)(z)$ satisfies the heat equation on $M$:
\begin{equation} {\pd\nu_f\over\pd t}(t, z) = A\nu_f(t, z), \qquad \nu_f(t=0,
z) = f(z). \tag {40}\end{equation}

\begin{lem} For every basic function $f$, we have 
$S_t f = T_t f$ for all $t \geq 0$. In particular, $S_t f$ is basic.
\end{lem}

\begin{proof} We have $\ddt S_t f = A S_t f$ in general.  Moreover, for {\it basic} $f$,

\begin{equation*}
\begin{split}
&{1\over2}\left( (\Delta_M + \kappa)f\right)\circ\pi =  (Af)\circ\pi = \wh{A} (f\circ\pi)\\
 &\qquad={1\over2}\sum_{k=1}^n Y_k^2 (f\circ\pi) = {1\over2}\sum_{k=p+1}^n Y_k^2 (f\circ\pi) = \wh{A^\bot}(f\circ\pi),\\
\end{split}
\end{equation*}
hence $\ddt T_t f = A T_t f$, where we have used the fact that $T_t f$ is basic for all $t$ (Lemma 3.4).  By uniqueness of solutions of the heat equation it follows that $S_t f = T_t f$.
\end{proof}

\begin{cor}The differential operator $A =
{1\over2}\Delta_M + {1\over2}\kappa$ leaves $C_{\text {b}}^\infty(M)$
invariant.
\end{cor}

\begin{proof} 
Recall that $\nu_f(t, z) = (S_tf)(z)$ and we have seen that $S_t$ preserves $C_{\text {b}}(M).$  Thus for $f \in C^\infty_{\text{b}}(M),$ each $\nu_f(t, \cdot)$ is basic and the
result follows by setting $t=0$ in (40).
\end{proof}
\bigbreak \par

By considering the scalarizations (\S 4), we can derive a result for $T_t$
acting on forms.

\begin{thm}  The infinitesimal generator of the transverse
semigroup $T_t$ acting on forms (32) is
\[ A = {1\over2}\Delta^\oplus,\]
where $\Delta^\oplus\theta = +\np_{e_i}(\np_{e_i}\theta) -
\np_{\np_{e_i}e_i}\theta,$ for any local orthonormal frame
$\{e_i\}$ in $\FO$ (summation on $i$ from $p+1$ to $n$ is understood).
In particular, $A$ preserves the basic complex.
\end{thm}
\begin{proof} The proof is analogous to that of the Cor.\ 5.3.
Equation (33) now holds componentwise for each function in the
scalarization $\{F_{\theta J}\}$ of $\theta.$  We need the fact that because
$Y_k$ is horizontal,
\begin{equation} Y_k F_{\theta\,J}(r) = (F_{\np\theta})_{J, k}(r).\tag {41} \end{equation}
This follows from a straightforward calculation, cf.\  Proposition
4.1 in [IW, Chap.\  V].  It also follows more conceptually from
the commutative diagram
\begin{equation*} \tag {42}
\begin{split}
\begin{CD}
  C^\infty(\FO, V^\Lambda)^G @> d+\rho^\Lambda_*(\omega_\bdot) >>\mathcal
A^1(\FO, V^\Lambda)_{basic} \\
     @V\alpha_0V\visom V  @V \alpha_1V\visom V \\
     \mathcal A^0(M, \mathcal A^j)  @> \np >> \mathcal A^1(M, \mathcal A^j)
\end{CD}
\end{split}
\end{equation*}
for the case of $j$-forms (see, e.g., [BGV, p.\  24]).  In (42)
$\frak g$ is the Lie algebra of the structure group $G
= O(p)\times O(q)$ of the principal bundle $\FO;$ $\frak g$ acts
by the differential $\rho^\Lambda_*$ of the representation
$\rho^\Lambda$ of $G$ on the vector space $V^\Lambda$ built up by
taking alternating tensor products of $\rho_0,$ the dual of the
standard representation of $G$ on $V = \Bbb R^p \oplus \Bbb R^q$
(recall the discussion around (8)); $C^\infty(\FO, V^\Lambda)^G$ is
the space of smooth $G$-equivariant maps; $\omega$ is the $\frak
g$-valued one-form (connection) corresponding to the covariant
derivative $\np.$  The scalarization $\{F_{\theta\,J}\}$ in (27)
gives the equivariant map in the upper left-hand corner of the
diagram, cf.\  (29).

For the second-order derivatives appearing in (33) (with the lower limit $k=1$ replaced by $k=p+1$), Eq.\  (41)
gives
\begin{equation} Y_k Y_k F_{\theta\,J}(r) = (F_{\np\np\theta})_{J, k,
k}(r).\tag {43}\end{equation}
From (32), (31), (43), and (33), with $\tnu$ replaced by
$\{F_{\theta_t\,J}\},$ it follows that
\[ {\pd\theta_t\over \pd t} = {1\over2}\Delta^\oplus\theta_t, \]
where $\theta_t \equiv T_t\theta$.

Arguing as in the proof of the above Corollary, but using this time Lemma
4.2, we see that $A$ preserves the basic complex.\qquad
\end{proof}
\bigbreak \par
\par

We close this section with a quick proof of the analog of Lemma 3.5 for forms.
\begin{lem} Let $\theta\in \Ab(M)$ be a basic $m$-form and let
$g, g'$ be two bundle-like metrics satisfying (3). Then
\[T_t\theta = T'_t\theta \text{  for all } t \geq 0.\]
\end{lem}
\begin{proof} We have from (32), (31), and the first equality in (27) that
$ T_t\theta (z) = \int_{\Omega_q} F_{\theta J}(R(t, r, \omega))P^W_0(d\omega) e_*^J(r) $
and
$ T'_t\theta (z) = \int_{\Omega_q} F_{\theta J}(R'(t, r', \omega))P^W_0(d\omega) {e'}_*^J(r').$
By Lemma 4.1(ii), only multi-indices $J$ with every component $>p$ appear in these equations.  We again choose $r'\in \FO'$ over $z\in M$ so that $C'(r') = C(r)$; thus $e_*^J(r) = {e'}_*^J(r')$.  Lemma 4.1(i) now permits us to repeat the proof of Lemma 3.5 with $f\circ\pi$ replaced by $F_{\theta J}$.
\end{proof}

Differentiating $T_t\theta = T'_t\theta$ at $t=0$, we obtain $ A\theta =
A'\theta \quad\text{for all basic forms }\theta, $ where $A, A'$ are
given by Theorem 5.4 for the metrics $g, g'$.
This result expresses a general invariance principle which
would be cumbersome to prove directly.

Finally, let us remark that the dependence on the homotopy class of $\gamma$ (i.e., covering-space phenomena associated with $\pi: \tL\to \mathcal L$) mentioned after Lemma 2.2 plays no role in this work.  For functions, this is because the projection $\pi$ appears in the definition (18) of $T_t$ and $S_t$; for basic forms $\theta$, it is because of Lemma 4.1(i).
\bigbreak \par
\section{THE FUNCTION $\phi$}
\bigbreak \par
\par
Because $P^W_0$ is a probability measure, the transition
semigroup $S_t$ (18) acts by contractions on $C(M),$ the Banach space of
continuous functions on $M$ with the sup norm.  The infinitesimal
generator $A = {1\over2}(\Delta_M +\kappa)$ acts on the smooth
functions $C^\infty(M) \subset C(M)$ and is closable.  The dual
semigroup $S_t^*$ acts on $C(M)^* = \text{Meas}(M),$
the Banach space of real-valued (signed) measures on $M$, and its infinitesimal
generator $A^*$ is a closed, densely defined operator on $C(M)^*$.
For $h\in C(M)$ smooth, $A^*h$ is given by the formal adjoint of $A$:
\begin{equation}A^*h = {1\over2}\left(\Delta_M\,h - \text{div} (h\kappa)\right)
= -\delta(dh - h\kappa)/2.\tag {44}\end{equation}  Here we regard $h$ as the
measure $h\,\text{dvol}_M$ on $M$, where $\text{dvol}_M$ is the
Riemannian volume element on $M$.
\par

Since we can work separately with each connected component,
there is no loss of generality in assuming $M$ to be connected as well as compact.  
It is then well known that the transition semigroup $S_t$ has a unique
invariant probability measure (see, e.g., [IW, Prop.\  V.4.5],
[Kun, Th.\  1.3.6], [N]), and by elliptic regularity this measure is of the
form $\phi \text{dvol}_g,$ with $\phi \geq 0$ smooth. We will need the
fact that $\phi>0$ everywhere.

\begin{prop} Let $M$ be compact and connected.  Then there exists a unique
probability measure $\mu(dz)$ invariant under $S_t$.  It is given
by $\phi\,\text{dvol}_M$, where $\phi \in C^\infty(M), \; \phi > 0
\text{ everywhere, and } A^* \phi = 0,$ i.e.,
\[ 0 = \delta(d\phi - \phi\kappa). \]
\end{prop}

\begin{proof} We sketch an argument [Ma]. Since $A$ is elliptic with vanishing zero-order part, its kernel reduces to the constants.
  By the index theorem, $\text{index}(A) = \text{index}(\Delta) = 0$, hence
$\text{dim ker}(A^*) = 1$. Choose $\phi \not\equiv 0$ with $A^*\phi = 0$; by elliptic theory, $\phi$ is smooth. 
 The associated  measure $\mu = \phi\, \text{dvol}_g$ on $M$ is invariant under the adjoint semigroup $S_t^*$, which like $S_t$ is a positivity-preserving contraction.  It then follows by a standard argument that we can take $\mu$ to be a positive measure, i.e., $\phi \geq 0$. 
 If $\phi$ were to vanish at some point $z_0\in M$, then writing out the equation 
$A^*\phi = 0$ in local coordinates and using the ellipticity of $A^*$, we see that all derivatives of $\phi$ through order two vanish at $z_0$.  Repeatedly differentiating the
 equation $A^*\phi = 0$, setting $z = z_0$, and proceeding by induction, we find that
 all derivatives of $\phi$ vanish at $z_0$.  Therefore, by Aronszajn's theorem $\phi
 \equiv 0$, a contradiction.  Alternatively, the results of [Bo] could also be used to show that $\phi > 0$.  
\end{proof}

\begin{dfn} Let $\psi > 0$ be smooth, $p = \text{ dim}\;\mathcal F$.
If $g'$ is obtained from $g$ by leaving $Q\equiv T\mathcal F^\bot$ unchanged
while rescaling $g$ along $T\mathcal F$ by $\psi^{2/p}$, so that $g' =
\psi^{2/p}g_{\mathcal F} \oplus g|_Q,$ we say that $g'$ is an $\mathcal F$-dilation
of $g$.
\end{dfn}

If $g$ is bundle-like (satisfies (3)), then clearly so is $g'$.
\par
Our immediate concern is with $\mathcal F$-dilations, for which we
will need to consider the long-time behavior $t\to \infty$.  Because the
generator 
$A = {1\over2}(\Delta_M + \kappa)$ of the transition semigroup
$S_t$ is not symmetric, we cannot argue as in the usual case of a
self-adjoint negative generator $A$, where $\lim_{t\to\infty}
e^{tA}\psi$ is the projection of the function or form $\psi$
onto its harmonic part.  But there is a substitute in the form of
the ergodic theorem ([Kun, Th.\  1.3.10]).  This holds for any
Feller semigroup $\{S_t\}$ for which the transition probability
$P_t(z, dw)$ is given by
\begin{equation}P_t(z, dw) = p_t(z, w)\text{vol}(dw)\tag {45}\end{equation}
for some strictly positive kernel $p_t(z, w)$ that is continuous
in $(t, z, w) \in (0, \infty) \times M^2.$  (We recall that the
transition probability $P_t(z, dw)$ is the measure defined by the
positive linear functional $f\mapsto S_tf(z),$ so that $S_tf(z) =
\int_M f(w) P_t(z, dw).$)
\par
The Feller condition is easily established (see, e.g., [Ma]). A proof that
the kernel $p(t, z, w) = p_t(z, w)$ exists and is continuous can be found
in [BGV, Th.\  2.23].  Since $S_tf(z) \geq 0$ for $f\geq 0$, we see 
that (45) holds with $p_t\geq 0$.  To show that $p_t > 0$, one can apply the strong
maximum principle; see, e.g., Theorem 3.1 in [Bo].  In fact, Bony's results hold quite generally for hypoelliptic operators and are thus more than we need here.  In
particular, strict positivity of the heat kernel for $T_t$ itself would follow if the latter were hypoelliptic, but this is hardly ever the case for Riemannian foliations.  So for technical reasons we work with $S_t$.  
\par
Thus the ergodic theorem applies to our situation and we conclude that for
any $f\in C(M)$ and $z\in M,$ \[ \lim_{t\to\infty} S_tf(z) = \int_M
f\phi\,\text{dvol}_g, \]
$\phi\,\text{dvol}_g$ being the unique invariant probability measure on
$M$ given by Proposition 6.1.
\par
We now dilate the bundle-like metric $g$ by $\phi$:
\begin{equation}g' = \phi^{2/p}g_{\mathcal F} \overset \bot \oplus g_Q. \tag {46}\end{equation}
Then
$ \text{dvol}_{g'} = \phi\,\text{dvol}_g.$

The new transition semigroup is $S'_t,$ and its infinitesimal generator $A'$
is given by $ A' = {1\over2}(\Delta_{g'} + \kappa')$,
where $\kappa' = \kappa - d_{1,0}\log\, \phi$ as follows from Rummler's formula (see for instance [Dom, Eq.\  (4.22)]).  By Lemmas 3.5 and 5.2, for all basic functions $f$
\begin{equation} S'_tf(z) = S_tf(z)\;\forall\,z\in M.\tag {47}\end{equation}
We note in passing that in the special case of dilations considered here it is not difficult to show directly that $A'f = Af$ for $f$ basic, hence $S'_tf = S_tf$ follows by the same uniqueness argument as in the proof of Lemma 5.2, thus avoiding Lemma 3.5.  However, Lemma 3.5 holds for arbitrary changes of metric subject to (3) and is useful in more general situations, as in Lemma 5.5.

Let us write $\phi'\,\text{dvol}_{g'}$ for the unique probability
measure on $M$ invariant under $S'_t;$ $\phi'$ is given by
Prop.\ 6.1.  For $f \in C_{\text{b}}(M)$ basic and $z\in M$ arbitrary, an
application of the ergodic theorem gives
\begin{align*}
     \lim_{t\to\infty}S'_tf(z) &= \int_M f\phi'\,\text{dvol}_{g'}
\cr
     &= \int_M f\phi'_{\text{b}'}\,\text{dvol}_{g'}\cr
     &= \int_M f\phi'_{\text{b}'} \phi\,\text{dvol}_g\cr
     &= \int_M f\phi'_{\text{b}'}\phib\,\text{dvol}_g,\cr
\end{align*}

and by (47) this is equal to
\begin{align*}
     \lim_{t\to\infty} S_tf(z) &= \int_M f\phi\,\text{dvol}_g\cr
     &= \int_M f\phib\,\text{dvol}_g.\cr
\end{align*}
Thus
\[ 0 = \int_M f[\phi'_{\text{b}'}\phib - \phib]\,\text{dvol}_g
\text{ for all basic }f,
\]
hence
\begin{equation} \phi'_{\text{b}'} \equiv 1, \tag {48}\end{equation}
since $\phib$ never vanishes [AL, Prop.\  2.2].

{\it Remark} 1. The above argument shows that
for any smooth basic function $\psi > 0$ on $M,$ there exists a
bundle-like metric $g',$ obtained from $g$ by a suitable $\mathcal
F$-dilation, such that $\psi = \phi'_{\text{b}'}.$

We recall that the exterior derivative $d$ preserves the basic
functions (and forms) $\Ab$.  Therefore, its adjoint $\delta$
preserves the $L^2$-orthogonal complement $\Abp.$  By Cor.\ 5.3, $A$
preserves the basic functions $C_{\text{b}}$, hence its adjoint $A^*$ leaves
$\Cbp$ invariant.  Writing $\phi = \phib + \phio$ as the sum of
its basic and orthogonal components, and using the fact that $\phib$ and $\phio$ are smooth,  we see that
$ A^*\phio \in \Cbp$.  Since $A^* f = -\delta(df - f\kappa)/2$ by
(44), we obtain
$ \delta(d\phio - \phio\kappa) \in \Cbp$.
Together with the argument leading to (48), this implies:
\bigbreak \par
\begin{thm} Let a bundle-like metric $g$
be given. Then there exists another bundle-like metric $g'$ on
$M$, obtained by a dilation of $g$ as in Eq.\  (46), with the
property that $\kb$ is basic-harmonic, i.e., $\delta_{\text{b}}\kb =
0 = d\kb.$
\end{thm}

\begin{proof}  By definition, $\delta_{\text{b}} = P_{\text{b}}
\circ\delta,$ where $P_{\text{b}}$ is the $L^2$ projection onto the
basic complex. According to  [AL, Cor.\ 3.5], $d\kb = 0.$  On
the other hand, using $A^*\phi = 0$ and $\phi =  \phib + \phio,$
we have
\[
 \delta(d\phib - \phib\kappa) = -\delta(d\phio - \phio\kappa) \in
\Cbp.
\]
Clearly, $\phib\ko \in \Abp,$ so $\delta(\phib\ko) \in \Cbp$ and therefore
\begin{equation} \delta(d\phib - \phib\kb) \in \Cbp. \tag {49}\end{equation}
Using the metric $g',$ we may suppose that $\phib$ is identically
equal to 1.  Then $\delta\kb \in \Cbp,$ i.e.,
$\delta_{\text{b}}\kb = 0.$
\end{proof}

{\it Remark} 2. This result is trivial if all basic
functions are locally constant, because any divergence
automatically integrates to zero.  In the contrary case, however, $\text{dim }d C_{\text{b}} = \infty$ and Theorem 6.2 solves an infinite-dimensional, global, nonlinear problem.

{\it Remark} 3. It is clear from Proposition 6.1 that
$\phi = \text{ const } \iff \delta\kappa = 0.$
Moreover,
\newcommand{\db}{\delta_{\text{b}}}
$ \phib = \text{ const } \iff \db\kappa = 0.$
The implication $\Rightarrow$ was shown in the proof of Theorem 6.2.
Conversely, suppose that $\db\kappa = 0.$ We always have
$ -\delta(d\phib - \phib\kb)\in \Cbp(M),$
but this is equal to
\begin{align*}
& \Delta\phib + \phib \delta\kb - \kb(\phib) \cr
&\qquad = \left(2A\phib - \kappa(\phib)\right) + \phib \delta\kb -
\kb(\phib)\cr
&\qquad = 2A\phib - 2\kb(\phib) - \ko(\phib) + \phib \delta\kb.\cr
\end{align*}
\newcommand{\Pb}{P_{\text{b}}}
The first two terms in the last line are in $\Cb(M),$ and by
hypothesis the last term is in $\Cbp(M).$ Moreover, $\Pb\,\ko(\phib) = 0,$
since $\Cbp\owns \delta(\phib\ko) = \phib \delta\ko - \ko(\phib)$ gives
$\Pb\ko(\phib) = \Pb(\phib\delta\ko) = \phib\,\Pb\delta\ko = 0.$ It follows that
$(A-\kb)\,\phib = 0$, hence by the maximum principle for elliptic operators,
$\phib = \text{ const}$.  

\par
Although the content of Theorem 6.2 is in no way changed, it takes a
somewhat nicer form ($\kb$ can be replaced by $\kappa$) if we assume the
truth of a long-standing conjecture asserting the existence of a bundle-like
metric with basic mean curvature. This conjecture has recently been proved
by \Dom.
\begin{cor} Let $M$ be a compact manifold
equipped with a Riemannian foliation, and let $g$ be a
bundle-like metric for which $\kappa$ is basic [Dom].  Then
$g$ can be dilated to obtain another bundle-like metric $g'$ for
which the mean curvature $\kappa'$ is basic-harmonic.
\end{cor}

\begin{proof}If $f$ is any smooth strictly positive function on $M,$ its
basic component is again smooth and strictly positive: $f_{\text{b}} > 0$ ([AL, Prop.\ 
2.2]).  Thus we need only dilate $g$ by $\phib;$ we saw in
(48) that $\phi'$ for the new metric $g'$ has constant basic
part. Since $\kappa' = \kappa - d_{1,0}\log\phib = \kappa -
d\log\phib$ is again basic, the result follows from the primed analog of (49), in which all quantities are for the metric $g'$.
\end{proof}

\bigbreak \par
The above corollary fits well with the Hodge decomposition for the basic
complex (see, e.g., [KT]). This gives an orthogonal decomposition
\[ \Ab(M) = \text{ im}\;d_{\text{b}} \oplus H_{\text{b}} \oplus
\text{ im}\;\delta_{\text{b}},\]
where $d_{\text{b}}$ is $d$ restricted to the basic forms and
$\delta_{\text{b}} = P_{\text{b}}\circ\delta,$ with $P_{\text{b}}$ the $L^2$
projection onto the basic complex.  The space $H_{\text{b}}$ consists of
those forms $\alpha$ satisfying $d_{\text{b}}\alpha = 0 = \delta_{\text{b}}\alpha$ and is
finite-dimensional.
Since $\kappa$ basic is equivalent to $d\kappa =0$, we know
{\it a priori} only that $\kappa\in \text{im}\;d_{\text{b}} \oplus
H_{\text{b}}$.  The Corollary  asserts that we can arrange for $\kappa$ to
lie in the finite-dimensional space $H_{\text{b}}$.  This result does
not seem to follow from the Hodge decomposition. For suppose that a
bundle-like metric $g$ with $\kappa$ basic has been found.  Then
$d\kappa = 0$ and we can write $\kappa = d_{\text{b}}f + h,$
where $f$ is basic and $h$ is basic-harmonic.  A natural thing to try is
to set $\lambda = e^f$ and dilate $g$ by $\lambda$ to get
$ \kappa' = \kappa - d_{1, 0} f = h.$
Then $\kappa'$ is again basic, but $h$ is in general not basic-harmonic for
the new metric $g'.$  More precisely, by Remark 3 and the argument leading
to (48), $h = \kappa'$ is basic-harmonic for $g'$ $\iff$ $\phi'_{\text {b}'}
= e^{-f} \phib$ is constant $\iff$ $\kappa = d_{\text{b}} \log\,\phib + h.$
\bigbreak \par
\section{AN EXAMPLE}
\bigbreak \par

We conclude with an example [Car].  Consider the manifold $M' = T\times
\Bbb R$ where $T$ is the 2-torus, and let $A\in SL(2, \Bbb Z)$ have trace
$> 2$.  Then $A$ has distinct real (irrational) eigenvalues $\lambda$ and
$1/\lambda$ with associated eigenvectors $V_1$ and $V_2$.  It defines an
orientation-preserving diffeomorphism of $T = \Bbb R^2/\Bbb Z^2$.  The
direction determined by $V_1$, say, defines a flow on $M$ by
\[ \psi_s( (x, y), t) = ( (x, y)+s V_1, t)\]
for $s\in \Bbb R$.  The integers $\Bbb Z$ act on $M'$ by $( (x, y), t)^m
= (A^m((x, y)), t+m),\;(x, y)$ a general point in $T$. Because $V_1$ is an
eigenvector of $A$, the flow defined by $\psi$ induces a one-dimensional
Riemannian foliation $\mathcal F$ on the compact quotient manifold
$M = M'/\Bbb Z$.
The nonconstant function $F([(x, y), t]) = \sin(2\pi t)$ is well-defined on
$M$ and is basic, hence the space $d_{\text{b}}(C_{\text{b}}(M))$ is
infinite-dimensional.  Carri\`ere shows that $(M, \mathcal F)$ admits a
transverse Lie structure modeled on the affine group $\Bbb R^2.$  This feature
enabled him to prove directly that the second basic cohomology group
vanishes: $H^2_{\text{b}} = 0.$
It follows that there exists no bundle-like metric for which
$\kappa = 0.$  For more details, we refer to Chapter 10 of [T].
Since $\kappa$ is nontrivial (in a rather strong sense) and nonconstant basic
functions exist, Theorem 6.2 has content in this case.
\par
Let us examine in more detail what our results say in the context of the above example.  
We take the leaf coordinate $x$ to be along $V_1$ and the transverse coordinates $y$
and $t$ to be along $V_2$ and the $t$ axis, respectively.  The local model space $\Bbb
R^2$ is identified with the affine group $GA(2)$ with group law $(y, t)\circ (y', t') = (\lambda^{-t} y'+y, t+t')$.
The transverse metric $g_T$ is taken to be any left-invariant metric on $GA(2)$.  This
amounts to assigning a metric arbitrarily at the identity element $(0, 0)$ and
transporting it by left multiplication.  Thus,
\begin{align*}
g_T\Big|_{(y, t)}(\lambda^{-t} {\pd\over\pd y}, \lambda^{-t} {\pd\over\pd y}) &= g_T\Big|_{(0, 0)}({\pd\over\pd y}, {\pd\over\pd y}), \cr
g_T\Big|_{(y, t)}(\lambda^{-t} {\pd\over\pd y}, {\pd\over\pd t}) &= g_T\Big|_{(0, 0)}({\pd\over\pd y}, {\pd\over\pd t}), \cr
g_T\Big|_{(y, t)}({\pd\over\pd t}, {\pd\over\pd t}) &= g_T\Big|_{(0, 0)}({\pd\over\pd t}, {\pd\over\pd t}).\cr
\end{align*}

In particular, there is no need to take ${\pd\over\pd y}$ and ${\pd\over\pd t}$ to be orthonormal at $(0, 0)$, though of course we could.  By construction, the metric $g_T$ is invariant under the
identification $(x, y, 0) = (\lambda x , \lambda^{-1} y, 1)\in T \times \Bbb R$ in the definition of $M$.  
\par
As stated after Eq.\  (3), given any Riemannian metric $g'$ on M, we obtain a bundle-like
metric compatible with $g_T$ by setting $g(X, Y) = g'(PX, PY) + g_T(\overline{X}, \overline{Y})$.  We could take $g'$ to come from the standard metric $g''$ on $T \times \Bbb R$, except within a buffer layer $T \times [1-c, 1)$, where $g''$ must be
deformed so as to be consistent with the identification $(x, y, 0) \sim (A(x, y), 1)$ and
give a well-defined metric $g'$ on the quotient $M$. Many other choices of $g'$ and hence
$g$ are possible; for instance, $T = S^1\times S^1$ and we could perturb the metrics on
each of the circle factors.  With the standard choice, ${\pd\over\pd y}$ will not be
orthogonal to ${\pd\over\pd x}$.  
\par
To find the mean curvature $\kappa$ in local coordinates we use the Koszul formula, which
requires computing the Lie brackets
$[e_1, e_2]$ and $[e_1, e_3]$ for an orthonormal frame $\{e_1, e_2, e_3\}$ with $e_1$ proportional to ${\pd\over\pd x}$ and $e_2$ and $e_3$ linear combinations of
${\pd\over\pd x}, {\pd\over\pd y}, \text{ and }{\pd\over \pd t}$, all coefficients
depending on the metric $g$.  This can be done explicitly but is not particularly illuminating.
Furthermore, there is little hope of actually finding the function $\phib(t)$ explicitly.

\newcommand{\dvol}{\text{dvol}_g}
\par
We now consider what the Corollary of Theorem 6.2 says in the present situation.  Since $\lambda$ is irrational, for each $t \in [0, 1)$ every leaf meeting the torus $T \times \{t\}$ is dense in it, hence the basic functions $F$ on $M$ depend only on
the $t$ coordinate and can be identified with the smooth functions on $\Bbb R^1$ with
period $1$. By [Dom, Theorem 4.18], given any $g_T$ there exists a bundle-like metric $g$
satisfying (3) for which $\kappa$ is basic.  Dilating by $\phib$, we can achieve in addition that
$\delta_{\text{b}}\kappa = 0$, i.e., $\int_M F'(t) (dt, \kappa) \text{dvol}_g = 0$ for
every smooth function $F$ with period 1 in $t$.  We set $h(t) = (dt, \kappa)$, which is a
basic function because $\kappa$ is basic and $g$ is bundle-like.  Taking $F(t)$ to be $\sin(2\pi mt)$ or
$\cos(2\pi mt)$ for $m\in \Bbb Z$, it follows that $\int_M \cos(2\pi mt) h(t)\dvol = 0$
and $\int_M \sin(2\pi mt) h(t)\dvol = 0$ for all $m$, except that $m=0$ must be excluded
in the first case.  Letting $F$ be any smooth periodic function with period 1 and
expanding $F$ in a Fourier series, we conclude that
\begin{equation} \int_M F(t) h(t)\dvol = C F_0, \tag {50}\end{equation}
where $C = \int_M h(t)\dvol $ and $F_0 = \int_0^1 F(t) dt$.  This equality extends by
continuity to periodic $F$ in $L^1[0, 1]$. 
\par
Replacing $dt$ by $-dt$ if necessary, we may
suppose that $C \geq 0.$  If $C = 0$ then (50) with $F = h$ shows that $h\equiv 0$, so let us   take $C \neq 0$.  Taking $F$ to be the characteristic function of $[\alpha, \beta]$, we find that $\int_{\alpha\leq t\leq \beta} h(t)\dvol = C(\beta - \alpha)$ for all $\alpha, \beta\in [0, 1]$.  It follows that 
\begin{equation} h(t)/C = {d\mu_L\over  d\mu}(t),  \tag {51}\end{equation}
the Radon--Nikodym derivative of Lebesgue measure on $[0, 1]$ with respect to the measure
$\mu$ defined on $[0, 1]$ by $\mu[\alpha, \beta] = \int_M \chi_{\{\alpha\leq t\leq\beta\}}(x, y, t)\;\dvol$.  Thus Cor.\ 6.3 is equivalent to the
assertion that $(dt, \kappa) = \int_M (dt, \kappa)\dvol\; d\mu_L/d\mu$.
\par
We observe parenthetically that unless $h\equiv 0$, we must have $h(t) > 0$ for all $t$,
since (50) and the monotone convergence theorem imply that $\text{Vol}(M) = C\,\int_0^1 {1\over h(t)}\,dt$.  Since $h$ is smooth, if it ever vanished then the integral could not
converge.  In particular, if $(dt, \kappa)$ ever vanishes (e.g., if $\kappa$ vanishes at
some point), then it vanishes identically.  We recall here Carri\`ere's result that there
exists no bundle-like metric for which $\kappa \equiv 0$.
\par
Passing to the general case, we expect Theorem 6.2 to be nontrivial for Riemannian
foliations of higher codimension.  Provided the maximum dimension of the leaf
closures is strictly less than the dimension of $M$, one expects
nonconstant basic functions to exist.

\bigbreak \par
\section{CONCLUDING REMARKS}
\bigbreak \par

Examination of the proof of Theorem 6.2 suggests that it might be possible to construct a proof that avoids probability theory.  Indeed, neither the ergodic theorem nor the existence of $\phi$ requires probability; moreover, Cor.\ 5.3 can be established independently (it can be deduced, for instance, from Proposition 4.1 in [PR]). 
 Furthermore, as noted after Eq.\  (47), no appeal to Lemma 3.5 is necessary.  However, the proof of (47), which is based on Lemma 5.2, does require Lemma 3.4 (and also the reduction to $\FO$).  Thus, as far as Theorem 6.2 is concerned, the role of the probability theory is confined to the proof that $S_t$ preserves the basic functions.  But this property is much stronger than Cor.\ 5.3. 
 Indeed, according to the theory of semigroups (see, e.g., [Y, Chap.\  IX]), it amounts to the following: For each $f\in C^2(M)$ and $\alpha \geq 0$, if $(1 - \alpha A) f$ is basic then $f$ is basic. 
 I do not see how to prove this without using Lemma 3.4. 
\par
Finally, it may be worth noting that there is a suggestive analogy between $\phi$ and the function $\lambda$ considered in [Dom], which satisfies $d_{1,0} \lambda - \lambda\ko \in \delta_{\mathcal F}{\mathcal A}^{1,1}$.  \Dom's proof might be greatly simplified, and its geometric content made more apparent, if $\lambda$ could be replaced by $\phi$.  
This was actually one of the original motivations for the present work.  One can show that $d_{1,0} \lambda - \lambda\ko \in \overline{\delta_{\mathcal F}{\mathcal A}^{1,1}}$, the Fr\'echet closure of the image $\delta_{\mathcal F}{\mathcal A}^{1,1}$, if and only if $\lambda_{\text{b}} = \text{const}$.  Hence if $\phi$ can replace $\lambda$ then we must have
$\phib = \text{const}$; that this can be achieved is the content of Theorem 6.2.  But we are unable to proceed further using our methods, because they give no control over the basic-orthogonal part $\phio$.   

\bigbreak \par
\centerline{ References} 
\bigbreak \par
\def\bb{\bigbreak}
\noindent [AL]\quad  J. A. Alvarez L\'opez, The Basic Component
of the Mean Curvature for Riemannian Foliations, {\it Ann. Global
Analysis Geom.} {\bf 10} (1992), 179-194.
\bb 
\noindent [BGV]\quad N. Berline, E. Getzler, and M. Vergne, {\it Heat
Kernels and Dirac Operators}, Springer-Verlag, 1992.
\bb
\noindent[Bi]\quad J.-M. Bismut, {\it M\'ecanique Al\'eatoire},
Lect.\  Notes Math.\  {\bf 866}, Springer-Verlag, 1981.
\bb
\noindent[Bo]\quad J.-M. Bony, Principe du maximum, in\'egalit\'e de Harnack et unicit\'e
du probl\`eme de Cauchy pour les op\'erateurs elliptiques d\'eg\'en\'er\'es, {\it Ann.\  Inst.\ 
Fourier Grenoble} {\bf 19} (1969), 277-304.
\bb
\noindent[Car]\quad Y. Carri\`ere, Flots Riemanniens, {\it Ast\'erisque} {\bf 116} (1984).
\bb
\noindent[Dom]\quad D. Dom\'{\i}nguez, Finiteness and Tenseness
Theorems for Riemannian Foliations, {\it Am.\  J.\  Math.\ } {\bf 120} (1998), 1237-1276.
\bb
\noindent[IW]\quad N. Ikeda and S. Watanabe, {\it Stochastic
Differential Equations and Diffusions}, 2nd ed., North Holland,
1989.
\bb
\noindent[KT]\quad F. Kamber and P. Tondeur, De-Rham--Hodge
Theory for Riemannian Foliations, {\it Math.\  Ann.\ } {\bf 277} (1987),
415-431.
\bb
\noindent[Kun]\quad H. Kunita, {\it Stochastic Flows and Stochastic
Differential Equations}, Cambridge Univ.\  Press, 1990.
\bb
\noindent[Ma]\quad A. Mason, An Application of Stochastic Flows to Riemannian Foliations, thesis, Univ.\  Ill.\  Urbana-Champaign (1997).
\bb
\noindent[Mo]\quad P. Molino, {\it Riemannian Foliations},
Birkh\"auser, Boston, 1988.
\bb
\noindent[N]\quad E. Nelson, The Adjoint Markoff Process, Duke
Math.\  J.\  {\bf 25} (1958), 671-690.
\bb
\noindent[PR]\quad E. Park and K. Richardson, The Basic Laplacian of a Riemannian Foliation, {\it Amer.\  J.\  Math.\ } {\bf 188} (1996), 1249-1275.
\bb
\noindent[T]\quad P. Tondeur, {\it Foliations on Riemannian Manifolds}, Springer-Verlag, 1988.
\bb
\noindent[Y]\quad K. Yosida, {\it Functional Analysis}, 6th edition, Springer-Verlag, 1980.
\bb
\end{document}